\documentclass[a4paper,11pt]{amsart}

\usepackage{amsmath,amsthm,amssymb,mathrsfs}

\usepackage[colorlinks, linkcolor=blue,anchorcolor=Periwinkle,
    citecolor=blue,urlcolor=Emerald]{hyperref}
\usepackage[all]{xy}

\usepackage{pifont}
\usepackage{graphicx}
\usepackage{psfrag}

\usepackage{setspace}

\usepackage{tikz-cd}

\usepackage{shorttoc}

\usepackage[normalem]{ulem}
\usepackage{tikz}
\usetikzlibrary{matrix,positioning,decorations.markings,arrows,decorations.pathmorphing,	
	backgrounds,fit,positioning,shapes.symbols,chains,shadings,fadings,calc}
\tikzset{->-/.style={decoration={  markings,  mark=at position #1 with
			{\arrow{>}}},postaction={decorate}}}
\tikzset{-<-/.style={decoration={  markings,  mark=at position #1 with
			{\arrow{<}}},postaction={decorate}}}

\numberwithin{equation}{section}

\newcommand{\ie}{{\em i.e.}\ }
\newcommand{\cf}{{\em cf.}\ }

\numberwithin{equation}{section}
\newtheorem{theorem}{Theorem}[section]
\newtheorem*{theorem*}{Theorem}

\newtheorem{lemma}[theorem]{Lemma}

\newtheorem{proposition}[theorem]{Proposition}
\newtheorem{corollary}[theorem]{Corollary}

\newtheorem*{conjecture*}{Conjecture}

\newtheorem{remark}[theorem]{Remark}
\newtheorem{definition}[theorem]{Definition}

\newcommand{\dimv}{\underline{\dim}\,}

\renewcommand{\mod}{\mathrm{mod}\,}

\newcommand{\Hom}{\mathrm{Hom}}

\newcommand{\opname}{\operatorname}

\newcommand{\Z}{\mathbb{Z}}

\newcommand{\cA}{{\mathcal A}}

\newcommand{\cM}{{\mathcal M}}
\newcommand{\cN }{{\mathcal N}}

\newcommand{\bA}{{\mathbb A}}
\newcommand{\bB}{{\mathbb B}}
\newcommand{\bC}{{\mathbb C}}
\newcommand{\bP}{{\mathbb P}}

\newcommand{\fC}{\textbf{C}}
\newcommand{\fM}{{\mathbf M}}
\newcommand{\fS}{{\mathbf S}}
\newcommand{\fT}{{\mathbf T}}
\newcommand{\fd}{{\mathbf d}}
\newcommand{\Arcseg}{\mathbf{Arcs}}
\newcommand{\Seg}{\mathbf{Seg}}
\newcommand{\ba}{\mathbf{a}}
\newcommand{\Int}{\mathbf{Int}}
\newcommand{\Intv}{\mathbf{\underline{Int}}}

\newcommand{\nn}{node[black]{$\bullet$}}

\renewcommand{\phi}{\varphi}

\renewcommand{\hat}[1]{\widehat{#1}}

\renewcommand{\tilde}[1]{\widetilde{#1}}

\newcommand{\vsim}{\vec{\sim}}
\newcommand{\nvsim}{\vec{\nsim}}

\newcommand{\kn}{\mathfrak{n}}
\newcommand{\km}{\mathfrak{m}}

\makeatletter
\DeclareRobustCommand{\cev}[1]{%
  {\mathpalette\do@cev{#1}}%
}
\newcommand{\do@cev}[2]{%
  \vbox{\offinterlineskip
    \sbox\z@{$\m@th#1 x$}%
    \ialign{##\cr
      \hidewidth\reflectbox{$\m@th#1\vec{}\mkern4mu$}\hidewidth\cr
      \noalign{\kern-\ht\z@}
      $\m@th#1#2$\cr
    }%
  }%
}
\makeatother

\begin{document}

\title[Intersection vectors]{Intersection vectors  over tilings with applications to gentle algebras and cluster algebras}\thanks{Partially supported by the National Natural Science Foundation of China (Grant No. 11971326, 12071315).}

\author[C. Fu]{Changjian Fu}
\address{Department of Mathematics, Sichuan University, Chengdu 610064, P.R.China}
\email{changjianfu@scu.edu.cn}
\author[S. Geng]{Shengfei Geng}
\address{Department of Mathematics, Sichuan University, Chengdu 610064, P.R.China}
\email{genshengfei@scu.edu.cn}

\subjclass[2000]{}

\begin{abstract}
It is proved that a multiset of permissible arcs over a tiling is uniquely determined by its intersection vector under a mild condition. This generalizes a classical result over marked surfaces with triangulations. We apply this result to study $\tau$-tilting theory of gentle algebras and denominator conjecture in cluster algebras.  In the case of gentle algebras, it is proved that different $\tau$-rigid $A$-modules over a gentle algebra $A$ have different dimension vectors if and only if $A$ has no even oriented cycle with full relations. For cluster algebras,  the denominator conjecture has been established for cluster algebras of type $\mathbb{A}\mathbb{B}\mathbb{C}$.

\end{abstract}

\dedicatory{Dedicated to Professor Bernhard Keller on the Occasion of his 60th Birthday.}

\maketitle

\tableofcontents

\section{Introduction and main results}

\subsection{Tilings}
The combinatorial geometry of marked surfaces has played an important role in the study of cluster algebras. In their seminal paper \cite{FST}, Fomin, Shapiro and Thurston initiated the study of a class of cluster algebras known as cluster algebras of surface type using the geometric models provided by marked surfaces with triangulations. Various concepts in cluster algebras can be interpreted geometrically in  marked surfaces. For instance, cluster variables have a one-to-one correspondence with permissible arcs, and clusters correspond to triangulations over the associated marked surface. This correspondence has allowed researchers to establish many conjectures and properties of cluster algebras of surface type.

 Due to the interaction between cluster algebras and representation theory of finite dimensional algebras, marked surfaces also  appear naturally  in the study of representation theory of finite dimensional algebras. This connection is explored in works such as \cite{BZ, L, OPS}. In particular, geometric models of certain (generalized) cluster categories or algebras have been constructed via marked surfaces \cite{BZ,OPS, BS,HZZ}. Similar to the case of cluster algebras, various objects in representation theory can be interpreted geometrically, which enable us to establish new results for the corresponding categories and algebras (cf.  \cite{APS} for instance).
 
 A tiling is a marked surface with a partial triangulation. Here and throughout, by a marked surface, we mean a marked surface without punctures. For a given tiling $(\fS,\fM,\fT)$, Baur and Sim\~{o}es \cite{BS} introduced a finite dimensional algebra $\Lambda^{\fT}$ called the {\it tiling algebra}. They  proved that tiling algebras are gentle algebras and every gentle algebra can be realized as a tiling algebra. Furthermore, a geometric model of module category of a gentle algebra has been constructed. In particular, every indecomposable module of a gentle algebra has a geometric interpretation via curves over the associated tiling.

 In this paper, we continue our investigation of $\tau$-tilting theory of gentle algebras initiated in \cite{FGLZ} via geometric models in the sense of \cite{BS,HZZ}. Here we focus on  $\tau$-rigid modules, not only the indecomposable ones. Hence we have to consider finite multisets of pairwise compatible permissible arcs.
 Recall that for a given tiling  $(\fS,\fM,\fT)$, the surface $\fS$ is divided by $\fT$ into a collection of tiles of type $\text{(I),(II),(III),(IV),(V)}$ (see Section \ref{ss:tiling} for the precise definition). Denote by $\mathscr{R}(\fT)$ the set of finite multisets of pairwise compatible permissible arcs over $(\fS,\fM,\fT)$. By construction, $\mathscr{R}(\fT)$ is in one-to-one correspondence with the set of isomorphism classes of $\tau$-rigid modules over the associated gentle algebra $\Lambda^\fT$.
  For a given multiset $\cM\in \mathscr{R}(\fT)$,  we have an integer vector $\Intv_{\fT}(\cM)$ called the {\it intersection vector} of $\cM$ with respect to $\fT$, which is a geometric interpretation of the dimension vector of the $\tau$-rigid module $M$ corresponding to $\cM$.
 
Our first main result  shows that a multiset $\cM$ is uniquely determined by its intersection vector under a mild condition.
 \begin{theorem}\label{t:main-result-1}
 	Suppose that the tiling $(\fS,\fM,\fT)$ has no tiles of type $(II)$ nor even-gons of type $(V)$.
 	For any $\cM, \cN\in \mathscr{R}(\fT)$ such that   $\Intv_{\fT}(\cM)=\Intv_{\fT}(\cN)$, we have $\cM=\cN$.
 \end{theorem}
 Theorem \ref{t:main-result-1} generalizes the following classical result for marked surfaces with triangulations: each arc is uniquely determined by its intersection vector with respect to the triangulation \cite{Mosher,T}.  Furthermore, it also strengthens \cite[Thoerem 1.1]{GyYu} for the corresponding setting. In order to prove Theorem \ref{t:main-result-1}, we obtain a local criterion to identify two finite multisets of permissible arcs, which is of independent interest.
 
 \subsection{Applications}
 Due to the connection between tilings and gentle algebras,  we immediately apply Theorem \ref{t:main-result-1} to study $\tau$-rigid modules over gentle algebras via the geometric model developed in \cite{BS}.
 
 \begin{theorem}\label{t:main-result-2}
Let $A$ be a finite dimensional gentle algebra over an algebraically closed field $k$. The following are equivalent:
\begin{itemize}
\item[(1)] Different $\tau$-rigid $A$-modules have different dimension vectors;
\item[(2)] The Gabriel  quiver $Q_A$ of $A$ does not admit an oriented cycle of even length with full relations.
\end{itemize}
\end{theorem}

As we have already mentioned, a special type of tilings and the associated gentle algebras provide categorifications for the corresponding cluster algebras of surface type. We then apply Theorem \ref{t:main-result-1} and Theorem \ref{t:main-result-2} to study the associated cluster algebras of surface type.  Let $\mathcal{A}$ be a cluster algebra with initial cluster variables $x_1,\dots,x_n$. Recall that a cluster monomial of $\mathcal{A}$ is a monomial of cluster variables belonging to the same cluster.  There are two kinds of  integer vectors called {\it denominator vectors} ($d$-vector for short) and {\it $f$-vectors}  associated to cluster monomials. We remark that the $d$-vector and $f$-vector of a cluster monomial depend on the choice of the initial cluster variables.  By definition, the $f$-vector of an initial cluster variables $x_i$ is the zero vector, while its $d$-vector is $-\mathbf{e}_i$, where $\mathbf{e}_1,\dots, \mathbf{e}_n$ is the standard basis of $\Z^n$.  In particular, for a cluster monomial involving initial cluster variables,  information provided by initial cluster variables is missing in its $f$-vector. In order to overcome this, we may  redefine the $f$-vector $\bar{\mathbf{f}}(x_i)$ of the initial cluster variable $x_i$ by setting $\bar{\mathbf{f}}(x_i)=-\mathbf{e}_i$. Then our first application of Theorem \ref{t:main-result-1}  in cluster algebras can be read as follows, which strengthens \cite[Corollary 1.5]{GyYu}  for the associated cluster algebras of surface type.

 \begin{theorem}\label{t:main-result-3}
 Let $(\fS,\fM)$ be  a marked surface with a triangulation $\fT$. Denote by $\mathcal{A}_\fT$ the associated cluster algebra. Different cluster monomials of $\mathcal{A}_{\fT}$ have different $f$-vectors. \end{theorem}

 Recall that a cluster algebra is of {\it finite type} if it has finitely many cluster variables. By Fomin and Zelevinsky's classification theorem \cite{FZ03}, cluster algebras of finite type are classified by finite root systems and are independent of the choice of coefficients. 
 According to \cite[Theorem 1.8]{Gyoda}, for a cluster algebra of finite type, the $f$-vector of a cluster monomial coincides with its $d$-vector.  On the other hand, cluster algebras of type $\mathbb{A}$ are surface type which can be realized by marked surfaces with triangulations.
 As a direct consequence of Theorem \ref{t:main-result-3}, we conclude that different cluster monomials have different $d$-vectors for a cluster algebra of type $\mathbb{A}$.  On the other hand, cluster algebras of type $\mathbb{C}$ can be categorified by  a class of gentle algebras \cite{FGL21a,FGL21b}. As an application of Theorem \ref{t:main-result-2} and a duality between cluster algebras of type $\mathbb{B}$ and $\mathbb{C}$, we prove a similar result for cluster algebras of type $\mathbb{B}$ and $\mathbb{C}$.
 In summary, we have the following.
\begin{theorem}\label{t:main-result-4}
Let $\mathcal{A}$ be a cluster algebra of type $\mathbb{A}$, $\mathbb{B}$ or $\mathbb{C}$, and $\mathbf{x}$ be an arbitrary cluster. Different cluster monomials have different  denominator vectors with respect to $\mathbf{x}$.
\end{theorem}

In particular, Theorem \ref{t:main-result-4} establishes the following {\bf denominator conjecture} for cluster algebras of type $\mathbb{A}\mathbb{B}\mathbb{C}$.
\begin{conjecture*}\cite[Conjecture 4.17]{FZ04} \cite[Conjecture 7.6]{FZ07} Let $\mathcal{A}$ be a cluster algebra and $\mathbf{x}$ be an arbitrary cluster. Different cluster monomials have different denominator vectors with respect to $\mathbf{x}$.
\end{conjecture*}
To our best knowledge, the denominator conjecture is still open widely. It has not yet been checked for cluster algebras of finite type. Indeed, the denominator conjecture has only been established for acyclic cluster algebras with respect to acyclic initial clusters. Note that not all clusters are acyclic for an acyclic cluster algebra of rank $>2$. Specifically, it has been verified for cluster algebras of rank $2$ by Sherman and Zelevinsky \cite{SZ}, for acyclic skew-symmetric cluster algebras with respect to an arbitrary acyclic cluster by Caldero and Keller \cite{CK06,CK08}, for cluster algebras of finite type with respect to a bipartite initial cluster by Fomin and Zelevinsky \cite{FZ07}.
By using the quantum cluster character, 
Rupel and Stella \cite{RS} further removed the condition of skew-symmetric type and verified the conjecture for all acyclic cluster algebras with respect to acyclic clusters. A weak version of the denominator conjecture has also been investigated, see \cite{GP, NS, FG, FGL21a} for instance.

The paper is structured as follows. In Section \ref{s:tiling}, we recall definitions related to tilings. We establish a local to  global criterion (Lemma \ref{l:key lemma}) in Section \ref{ss:local-to-global}. As an application,  we prove Theorem \ref{t:main-result-1} in Section \ref{s:intersection-vector}. Section \ref{s:dimension-vector} is devoted to proving Theorem \ref{t:main-result-2} by using Theorem \ref{t:main-result-1} and the geometric interpretation of $\tau$-tilting theory of gentle algebras. In Section \ref{s:app-cluster}, we apply Theorem \ref{t:main-result-1} and Theorem \ref{t:main-result-2} to investigate corresponding cluster algebras.  After recalling basic results of cluster algebras and cluster algebras of surface type, we present the proof of Theorem \ref{t:main-result-3} in Section \ref{ss:cluster-algebra-surface-type}.  The proof of Theorem \ref{t:main-result-4} is given in Section \ref{ss:proof-4}.

\subsection*{Acknowledgments}
We are very grateful to Prof. Yu Zhou for helpful discussion. 

\section{Tilings and permissible arcs}\label{s:tiling}

\subsection{Tilings}\label{ss:tiling}

A {\it marked surface} is a pair $(\fS, \fM)$, where $\fS$ is a connected oriented Riemann surface with non-empty boundary $\partial \fS$, and $\fM\subset \partial \fS$ is a finite set of marked points on the boundary. We assume that, if $\fS$ is a disc, then $\fM$ has at least four marked points.
A connected component of $\partial \fS$ is a {\it boundary component} of $\fS$.  A boundary component $B$ of $\fS$ is {\it unmarked} if $\fM\cap B=\emptyset$. We allow unmarked boundary components of $\fS$. A {\it boundary segment} is the closure of a component of $\partial \fS\backslash \fM$.

An {\it arc} on $(\fS, \fM)$ is a continuous map $\gamma: [0,1]\to \fS$ such that
\begin{itemize}
    \item[$\circ$] $\gamma(0),\gamma(1)\in \fM$ and $\gamma(t)\in \fS\backslash \fM$ for $0<t<1$;
    \item[$\circ$] $\gamma$ is neither null-homotopic nor homotopic to a boundary segment.
    
\end{itemize}
We always consider arcs on $\fS$ up to homotopy relative to their endpoints.
The inverse of an arc $\gamma$ on $\fS$ is defined as $\gamma^{-1}(t)=\gamma(1-t)$ for $t\in [0,1]$. An arc $\gamma$ is called a {\it loop} if $\gamma(0)=\gamma(1)$.
Denote by $\fC(\fS)$ the set of equivalence classes of arcs on $\fS$ under the equivalence relation given by taking inverse. 

For any arcs $\gamma_1, \gamma_2\in \fC(\fS)$, whenever we consider their intersections, we always assume that they are representatives in their homotopy classes such that the number of their intersections is minimal.
The {\it intersection number} between $\gamma_1$ and $\gamma_2$ is defined to be
\[\Int(\gamma_1,\gamma_2):=|\{(t_1,t_2)~|~0<t_1,t_2<1, \gamma_1(t_1)=\gamma_2(t_2)\}|.
\]
Clearly, $\Int(\gamma_1,\gamma_2)=\Int(\gamma_2,\gamma_1)$.  We say that  $\gamma_1$ and $\gamma_2$ are {\it compatible} if $\Int(\gamma_1,\gamma_2)=0$.  An arc $\gamma\in \fC(\fS)$ is said to be without self-intersections if $\Int(\gamma,\gamma)=0$.

 Since we are only interested in arcs without self-intersections. For simplicity of terminology, we make the following convention throughout the paper.

 \noindent{\bf Convention A.}
 By an arc, we always mean an arc without self-intersections.

\begin{definition}
 A {\it  partial triangulation} $\fT$ of $(\fS, \fM)$ is a  collection of pairwise compatible arcs in $\fC(\fS)$. In this case, the triple $(\fS, \fM,\fT)$ is called a {\it tiling}.
\end{definition}

\begin{figure}[h]
\begin{minipage}[t]{0.15\linewidth} 
\begin{tikzpicture}[xscale=0.8,yscale=0.8]

\draw[thick,fill=black!20] (2,0.17)arc (250:290:3);

\fill(3,0) circle(2pt);

\draw[thick,color=blue] (3,0) .. controls (2,-1) and (2,-2) .. (3,-2);
\draw[thick,color=blue] (3,0) .. controls (4,-1) and (4,-2) .. (3,-2);

\draw[thick,fill=black!20] (3,-1.3) circle(8pt);
\node at (3,-3){Type (I)};

\end{tikzpicture}
\end{minipage}
\begin{minipage}[t]{0.16\linewidth} 
\centering
		\begin{tikzpicture}[xscale=0.8,yscale=0.8]
\draw[thick,fill=black!20] (2,0.17)arc (250:290:3);

\fill(3,0) circle(2pt);

\draw[thick,fill=black!20] (2,-2.17)arc (110:70:3);
\fill(3,-2) circle(2pt);
\draw[thick,fill=black!20] (3,-1) circle(8pt);
\draw[thick,color=blue] (3,0) .. controls (2,-1) and (2.5,-1.5) .. (3,-2);
\draw[thick,color=blue] (3,0) .. controls (4.5,-1) and (4,-1.5) .. (3,-2);
\node at (3,-3){Type (II)};

\end{tikzpicture}
\end{minipage}
\begin{minipage}[t]{0.15\linewidth} 
\centering
		\begin{tikzpicture}[xscale=0.8,yscale=0.8]
\draw[thick,fill=black!20] (2,0.17)arc (200:240:3);
\draw[thick,fill=black!20] (2,0.17)arc (340:300:3);
\draw[blue,thick] (0.7,-1.4) to (3.3,-1.4);	
\fill(2,0.17) circle(2pt);
\fill(0.7,-1.4) circle(2pt);
\fill(3.3,-1.4) circle(2pt);
\node at (2,-3){Type (III)};

\end{tikzpicture}
\end{minipage}
\begin{minipage}[t]{0.25\linewidth} 
\centering
		\begin{tikzpicture}[xscale=0.5,yscale=0.5]
					
	 \draw[thick,fill=black!40] (-1.3,-2.15)arc (120:60:1.5);
			
	\draw[blue,thick] (0,2) to (-1,2);
	\draw[blue,thick] (-1,2)\nn to (-2,1.5)\nn;
	\draw[blue,thick] (-2,1.5)\nn to (-2.7,0.5)\nn;
			
	\draw[blue,thick] (-1,-2)\nn to (-2,-1.5)\nn;
	\draw[blue,thick] (-2,-1.5)\nn to (-2.7,-0.5)\nn;
		
	\draw[blue,thick] (0,2)\nn to (1,1.5)\nn;
	\draw[blue,thick] (1,1.5)\nn to (1.7,0.5)\nn;
			
	\draw[blue,thick] (0,-2)\nn to (1,-1.5)\nn;
	\draw[blue,thick] (1,-1.5)\nn to (1.7,-0.5)\nn;
        \draw[blue](-2.7,-0.4)node[above]{$\vdots$}(1.7,-0.4)node[above]{$\vdots$};

	\node at (0,-3){ Type (IV)};
			\end{tikzpicture}
		
		\end{minipage}%
\begin{minipage}[t]{0.2\linewidth} 
\centering
		\begin{tikzpicture}[xscale=0.5,yscale=0.5]

		\draw[blue,thick] (0,2) to (-1,2);
		\draw[blue,thick] (-1,2)\nn to (-2,1.5)\nn;
		\draw[blue,thick] (-2,1.5)\nn to (-2.7,0.5)\nn;
	        \draw[blue,thick] (0,-2)\nn to (-1,-2)\nn;
		\draw[blue,thick] (-1,-2)\nn to (-2,-1.5)\nn;
		\draw[blue,thick] (-2,-1.5)\nn to (-2.7,-0.5)\nn;
		
		\draw[blue,thick] (0,2)\nn to (1,1.5)\nn;
		\draw[blue,thick] (1,1.5)\nn to (1.7,0.5)\nn;
			
		\draw[blue,thick] (0,-2)\nn to (1,-1.5)\nn;
		\draw[blue,thick] (1,-1.5)\nn to (1.7,-0.5)\nn;

		\draw[blue](-2.7,-0.4)node[above]{$\vdots$}(1.7,-0.4)node[above]{$\vdots$};

	\node at (0,-3){ Type (V)};
			\end{tikzpicture}

	\end{minipage}	
		
\caption{Basic tiles}\label{f:tile-type 2}		
		
	\end{figure}
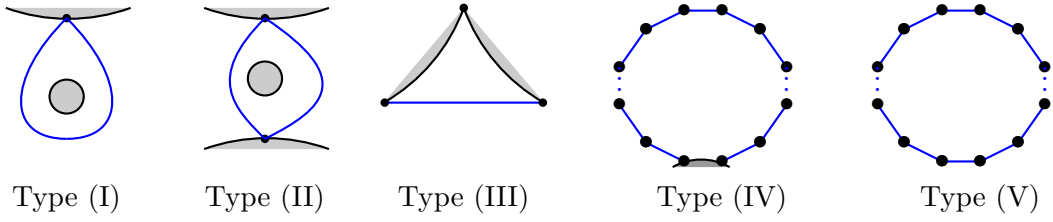

Let $(\fS, \fM, \fT)$ be a tiling.  Then $\fS$ is divided by $\fT$ into a collection of regions  (also
called {\it tiles}). We consider the following types of tiles (see Figure~\ref{f:tile-type 2}):
 \begin{enumerate}
     \item[(I)] monogons, that is, loops, with exactly one unmarked boundary component in their interior;
     \item[(II)] digons with exactly one unmarked boundary component in their interior;
  \item[(III)] three-gons bounded by two boundaries and one arc in $\fT$ and whose interior contains no unmarked boundary component of $\fS$;
       \item[(IV)] $m$-gons, with $m\ge 3$ vertices and whose edges are arcs in $\fT$ and one  boundary segment, and whose interior contains no unmarked boundary component of $\fS$;
    \item[(V)] $m$-gons, with $m\ge 3$ vertices and whose edges are arcs in $\fT$ and whose interior contains no unmarked boundary component of $\fS$.
 \end{enumerate}
 Throughout this paper, we also make the following convention.
 
  \noindent{\bf Convention B.}
 By a tiling $(\fS,\fM,\fT)$, we mean a tiling with tiles of type $\text{(I)}$-$\text{(V)}$.

 \subsection{Permissible arc} Let $(\fS,\fM,\fT)$ be a tiling and $\Delta$ a tile.
An {\it irreducible arc segment} in $\Delta$ is a curve $\eta: [0,1]\to \Delta$ such that  $\eta(0)$ and $\eta(1)$ are in the edges of $\Delta$
and $\eta(t)(0<t<1)$ are  in the interior of $\Delta$. 
 An irreducible arc segment $\eta$ in $\Delta$ is called {\it permissible} with respect to $\fT$ if it  satisfies one of the following conditions (cf. \cite[Definition 2.1]{HZZ} and compare \cite[Definition~3.1]{BS}).
	\begin{enumerate}
		\item[(P1)] One endpoint $P$ of $\eta$ is in $\fM$ and the other $Q$ is in the interior of a non-boundary edge, say $\gamma$ of $\Delta$, such that $\eta$ is not isotopic to a segment of an edge of $\Delta$ relative to their endpoints, and after moving $P$ along the edges of $\Delta$ in anti-clockwise order to the next marked point, say $P'$, the new arc segment obtained from $\eta$ is isotopic to a segment of $\gamma$ relative their endpoints. See Figure~\ref{f:P1} for all the possible cases of permissible irreducible arc segments satisfying (P1).
		
		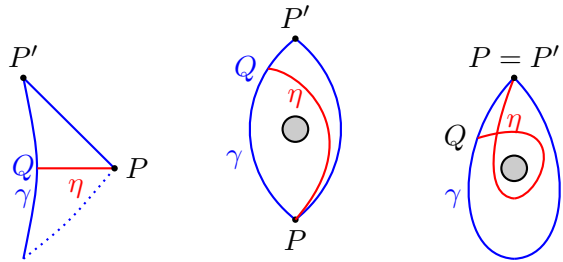
\begin{figure}[h]
			\begin{tikzpicture}[xscale=0.8,yscale=0.8]
				\draw[thick,color=blue] (1,0) .. controls (1.3,1.5) and (1.3,1.5) .. (1,3);
				\draw[thick,color=blue] (1,3)node[black,above]{$P'$} .. controls (2,2) and (2, 2) .. (2.5,1.5);
				\draw[thick,color=blue,dotted] (1,0) .. controls (1.75, 0.6) and (1.75,0.6) .. (2.5,1.5);
				\draw[thick,color=red] (2.5,1.5)node[black,right]{$P$} tonode[below]{$\eta$} (1.23,1.5);
				\fill(2.5,1.5) circle(1.5pt); \fill(1,3) circle(1.5pt);
				\draw[blue] (1,1)node{$\gamma$} (1,1.5)node[blue]{$Q$};	
			\end{tikzpicture}\qquad
			\begin{tikzpicture}[xscale=0.8,yscale=0.8]
				\draw[thick,color=blue] (7,0) .. controls (8,0.8) and (8,2.2) .. (7,3)node[black,above]{$P'$};	
				\draw[thick,color=blue] (7,0) .. controls (6,0.8) and (6,2.2) .. (7,3);	
				\fill(7,0) circle(1.5pt);		\fill(7,3) circle(1.5pt);	
				\draw[thick,fill=black!20] (7,1.5) circle(6pt);
				\draw[thick,color=red] (7,0)node[black,below]{$P$} .. controls (8.3,1.5) and (7,2.5) .. (6.55,2.5)node[blue,left]{$Q$};
				\draw[red] (7,2)node{$\eta$};
				\draw[blue] (6,1)node{$\gamma$};
			\end{tikzpicture}\qquad
		\begin{tikzpicture}[xscale=0.8,yscale=0.8]
				\draw[thick,color=blue] (13,0) .. controls (14,0.05) and (14,2) .. (13,3)node[black,above]{$P=P'$};
				\draw[thick,color=blue] (13,0) .. controls (12,0.05) and (12,2) .. (13,3);
				\draw[thick,color=red] (13,1) .. controls (12.5,1.05) and (12.6,2) .. (13,3);	
				\draw[thick,color=red] (13,1) .. controls (13.5,1.05) and (14,2.5) .. (12.4,2)node[black,left]{$Q$};	
				\draw[thick,fill=black!20] (13,1.5) circle(6pt);
				\fill(13,3) circle(1.5pt);
				\draw[red] (13,2.25)node{$\eta$};
				\draw[blue] (12,1)node{$\gamma$};
			\end{tikzpicture}
			\caption{Condition (P1)}\label{f:P1}
		\end{figure}
		
		\item[(P2)] The endpoints of $\eta$ are in the interiors of non-boundary edges $x,y$ (which are possibly not distinct) of $\Delta$ such that $\eta$ has no self-intersections, $x,y$ have a common endpoint $p_\eta\in \fM$ and $\eta$ cuts out an angle from $\Delta$ as shown in Figure \ref{f:PAS}. We denote by $\triangle(\eta)$ the local triangle cutting out by $\eta$ and $\angle \eta$ the angle opposite to $\eta$ in the local triangle $\triangle(\eta)$.
		\begin{figure}[htpb]
			\begin{tikzpicture}[xscale=0.8,yscale=0.8]
				\draw[thick,fill=black!20] (9,0.17)arc (250:290:3);
				\node at (10,0.2){$p_{\eta}$};
				\fill(10,0) circle(1.5pt);
				\draw[thick,color=blue] (10,0) .. controls (8.6,-0.9) and (8.4,-1.1) .. (8,-2.2);
				\draw[thick,color=blue,dashed] (8,-2.2) -- (7.8,-3);
				\draw[thick,color=blue] (10,0) .. controls (10.4,-0.9) and (10.6,-1.1) .. (11,-2.2);
				\draw[thick,color=blue,dashed] (11,-2.2) arc (20:-40:1);
				\draw[thick,color=red] (8,-1) .. controls (9,-1.5) and (10,-1) .. (12,-1);
				\draw[thick,color=red,dashed] (12,-1) arc (90:50:1.5);
				\draw[thick,color=red,dashed] (8,-1) -- (7,-0.5);
				\node at (9.5, -1.5){$\eta$};
				\node at (7.5,-3){{\color{blue}$x$}};
				\node at (10.5, -3){{\color{blue}$y$}};
				\draw (9.7,-.8)node{$\triangle(\eta)$};
			\end{tikzpicture}
			\caption{Condition (P2)}\label{f:PAS}
		\end{figure}
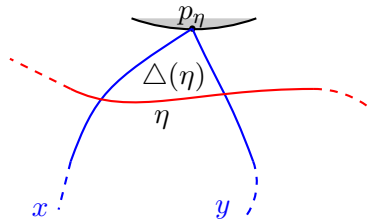
	\end{enumerate}
	
	For any arc $\gamma\in\mathbf{C}(\fS)$, we always assume that $\gamma$ is in a minimal position with $\fT$. The arc $\gamma$ is divided into irreducible arc segments by $\fT$.
	\begin{definition}[{\cite[Definition 3.1]{BS} and \cite[Definition 2.6]{HZZ}}]
		An arc $\gamma$ on $\fS$ is called {\it permissible} (with respect to $\fT$) if each irreducible arc segment of $\gamma$ divided by $\fT$ is permissible.
	\end{definition}
 
 \begin{remark}
 By our convention on arcs, the third permissible arc segment in Figure~\ref{f:P1} cannot be an arc segment of a permissible arc.
 \end{remark}

 We denote by $\bP(\fT)\subset\mathbf{C}(\fS)$ the set of permissible arcs on $\fS$ and by $\mathscr{R}(\fT)$ the set of  finite multisets of pairwise compatible permissible arcs.

\subsection{Equivalence of permissible arc segments}
Let $(\fS,\fM,\fT)$ be a tiling.
Let $\gamma$ be a permissible arc with a fixed orientation which is
divided into irreducible arc segments $\gamma_1,\dots, \gamma_m$ in order. For $1\leq i_1\leq i_1+k\leq m$, the concatenation $\gamma^{\langle k\rangle}:=\gamma_{i_1}\circ \gamma_{i_1+1}\circ\cdots\circ \gamma_{i_1+k}$ of $\gamma_{i_1}, \dots, \gamma_{i_1+k}$ is called  a (permissible) {\it arc segment of length $k$} for  $\gamma$. In particular, arc segments of length one are exactly the irreducible arc segments.
The orientation of $\gamma$ induces an orientation
for any arc segment $\gamma^{\langle k\rangle }$ of length $k$ of $\gamma$. Conversely, if $\gamma^{\langle k\rangle}$ is an arc segment of length $k$ for $\gamma$, then any orientation of $\gamma^{\langle k\rangle}$ can be extended
to an orientation of  $\gamma$ and induces an orientation for each irreducible arc segments of $\gamma^{\langle k\rangle}$. 

\noindent{\bf Convention C.}
In the following, when we write an arc segment as $\delta=\delta_1\circ\cdots\delta_l$, we make the convention that $\delta_i(1)=\delta_{i+1}(0)$ ($1\leq i<l$) and the orientation of $\delta$ is chosen such that $\delta(0)=\delta_1(0)$ and $\delta(1)=\delta_l(1)$.

Let $\Delta$ be a tile of $(\fS,\fM,\fT)$.
Let $\eta$ and $\gamma$ be two irreducible permissible arc segments in $\Delta$, we say that $\eta$ and $\gamma$ are {\it special homotopic}, denote by $\eta\ \vec{\sim}\ \gamma$, if $\eta$ is homotopic to $\gamma$ and we are in one of the following situations:
\begin{itemize}
   
    \item[(sh1)] $\eta(0)=\gamma(0)$ (resp. $\eta(1)=\gamma(1)$) belongs to $\fM$ and $\eta(1),\gamma(1)$ (resp. $\eta(0), \gamma(0)$) belong to the interior of an edge of $\Delta$ (cf. Condition (P1));
    \item[(sh2)] $\eta(0)$ and $\gamma(0)$ belong  to  the interior of an edge of $\Delta$, while $\eta(1)$ and $\gamma(1)$ are in  the interior of an edge of $\Delta$ (cf. Condition (P2)).
\end{itemize}
We say that $\eta$ is {\it equivalent} to $\gamma$ and denote is as $\eta\sim \gamma$, if $\eta\ \vec{\sim}\ \gamma$  or $\eta\ \vec{\sim}\ \gamma^{-1}$. 
We extend the equivalence relation of irreducible permissible arc segments to all the arc segments of permissible arcs. Let $\eta$ and $\gamma$ be  arc segments of length $l$ for some permissible arcs respectively, and denote them as $\eta=\eta_1\circ\cdots\circ \eta_l$ and $ \gamma=\gamma_1\circ\cdots\circ \gamma_l$, where $\eta_i$ and $\gamma_i$ are irreducible permissible arc segments. We say that $\eta$ is {\it equivalent} to $\gamma$ and denote by $\eta\sim \gamma$ if either $\eta_i\ \vec{\sim}\ \gamma_i(1\leq i\leq l)$ or $\eta_i\ \vec{\sim}\ \gamma_{l+1-i}^{-1}(1\leq i\leq l)$. We denote $\eta\ \vec{\sim}\  \gamma$ to indicate the case that $\eta_i\ \vec{\sim}\ \gamma_i$ ($1\leq i\leq l$) and $\eta\  \cev{\sim}\ \gamma$ to the other case.
In the following, we always consider  arc segments up to equivalence.
The definition of {\it compatible} of arcs can be extended to arc segments in an obvious way.

\begin{figure}[h]
\begin{tikzpicture}[xscale=0.9,yscale=0.9]

\draw[thick,fill=black!20] (2,0.17)arc (250:290:3);

\fill(3,0) circle(2pt);

\draw[thick,color=black] (3,0) .. controls (2,-1) and (2,-2.3) .. (3,-2.3);
\draw[thick,color=black] (3,0) .. controls (4,-1) and (4,-2.3) .. (3,-2.3);

\draw[thick,fill=black!20] (3,-1.6) circle(6pt);

\draw[thick,fill=black!20] (6,0.17)arc (250:290:3);

\fill(7,0) circle(2pt);

\draw[thick,color=black] (7,0) .. controls (6,-1) and (6,-2.3) .. (7,-2.3);
\draw[thick,color=black] (7,0) .. controls (8,-1) and (8,-2.3) .. (7,-2.3);
\draw[thick,fill=black!20] (7,-1.6) circle(6pt);

\draw[thick,red] (3.5,-0.7)  to (7.5,-0.7);

	\draw[red] (5.3,-0.8)node[above]{\tiny $\gamma^{(1)}$} (7,-0.8)node[above]{\tiny $\gamma^{(2)}$}; 

\draw[thick,blue] (2.4,-0.9)  to (6.4,-0.9);
\draw[blue] (5.3,-0.8)node[below]{\tiny $\delta^{(2)}$} (3,-0.8)node[below]{\tiny $\delta^{(1)}$}; 

\draw[black] (3,0)node[above]{\tiny $A$} (7,0)node[above]{\tiny $B$}; 
\end{tikzpicture}
\caption{Non-equivalence of arc segments}\label{f:nsim-loop}	
\end{figure}
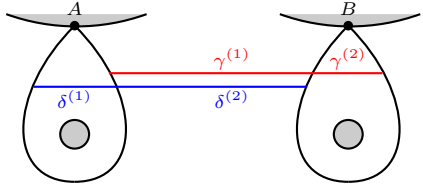

The following statements are direct consequences of the definition of equivalence of arc segments.
\begin{lemma}\label{l:nsim}
    Let $\gamma$ be an arc segments of length $l$ of a permissible arc. Then $\gamma\nvsim \gamma^{-1}$.
\end{lemma}
\begin{lemma}\label{l:nsim-loop}
   Let $\gamma=\gamma^{(1)}\circ \gamma^{(2)}$ and $\delta=\delta^{(1)}\circ\delta^{(2)}$ be arc segments of length $l$ for certain permissible arc as in Figure \ref{f:nsim-loop}  for some $l\ge 2$, where $\gamma^{(2)},\delta^{(1)}$ are of length $1$ and $\gamma^{(1)},\delta^{(2)}$ are of length $l-1$ such that $\gamma^{(1)}\vsim \delta^{(2)}$. Then $\gamma\nsim \delta$.
\end{lemma}

\begin{lemma}\label{l:equiv}
Let $\alpha=\alpha^{(1)}\circ \alpha^{(2)}\circ \alpha^{(3)}$ and $\beta=\beta^{(1)}\circ\beta^{(2)}\circ \beta^{(3)}$ be arc segments of length $l+1$ for certain permissible arc as in Figure \ref{f:equivalence-arc-seg}, where $\alpha^{(1)},\alpha^{(3)},\beta^{(1)},\beta^{(3)}$ are of length $1$ and $\alpha^{(2)},\beta^{(2)}$ are of length $l-1$. Suppose that $\alpha^{(2)}\vsim \beta^{(2)}$, then 
\begin{itemize}
    \item[(1)] $\alpha^{(2)}\circ \alpha^{(3)}\nsim \beta^{(1)}\circ\beta^{(2)}$;
    \item[(2)] $\alpha^{(2)}\circ \alpha^{(3)}\sim \beta^{(2)}\circ\beta^{(3)}$ if and only if $\alpha^{(2)}\circ \alpha^{(3)}\vsim \beta^{(2)}\circ\beta^{(3)}$
\end{itemize}
    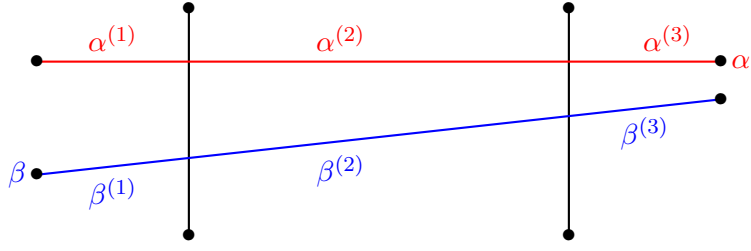
\begin{figure}[htpb]
		\begin{tikzpicture}
			

\draw[black,thick] (-3,0)\nn to (-3,-3)\nn;
\draw[black,thick] (2,0)\nn to (2,-3)\nn;


     \draw[red,thick] (-5,-0.7)\nn to (4,-0.7)\nn;
						\draw[red] (4,-0.7)node[right]{$\alpha$};	
	\draw[red] (-4,-0.7)node[above]{$\alpha^{(1)}$} (-1,-0.7)node[above]{$\alpha^{(2)}$}(3.3,-0.7)node[above]{$\alpha^{(3)}$};

 \draw[blue,thick] (4,-1.2)\nn to (-5,-2.2)\nn;
			\draw[blue] (-5,-2.2)node[left]{$\beta$};
	\draw[blue] (-4,-2.1)node[below]{$\beta^{(1)}$} (-1,-2.5)node[above]{$\beta^{(2)}$}(3,-1.3)node[below]{$\beta^{(3)}$};

		\end{tikzpicture}
		\caption{Equivalence of arc segments}\label{f:equivalence-arc-seg}
	\end{figure}
\end{lemma}
\begin{proof}
Let us first consider the case $l=1$.
In this case, $\alpha^{(3)}$ and $\beta^{(1)}$ belong to different tiles, which implies that $\alpha^{(3)}\nsim \beta^{(1)}$. On the other hand, $\alpha^{(3)}\sim \beta^{(3)}$ if and only if $\alpha^{(3)}\vsim \beta^{(3)}$ by Lemma \ref{l:nsim}.

Now assume that $l\geq 2$. If $\alpha^{(2)}\circ \alpha^{(3)}\sim \beta^{(1)}\circ\beta^{(2)}$, we  have either $\alpha^{(2)}\circ \alpha^{(3)}\vsim \beta^{(1)}\circ\beta^{(2)}$ or $\alpha^{(2)}\circ \alpha^{(3)}\cev{\sim} \beta^{(1)}\circ\beta^{(2)}$. If $\alpha^{(2)}\circ \alpha^{(3)}\vsim \beta^{(1)}\circ\beta^{(2)}$, we are in the situation of Lemma \ref{l:nsim-loop}, a contradiction. If $\alpha^{(2)}\circ \alpha^{(3)}\cev{\sim} \beta^{(1)}\circ\beta^{(2)}$, then $\alpha^{(2)}\cev{\sim}\beta^{(2)}$. By the assumption, we obtain $\alpha^{(2)}\vsim {\alpha^{(2)}}^{-1}$, which contradicts to Lemma \ref{l:nsim}. This finishes the proof of $(1)$.

For $(2)$, only the implication ``$\Rightarrow$" requires a proof. Suppose that $\alpha^{(2)}\circ \alpha^{(3)}\vsim {\beta^{(3)}}^{-1}\circ {\beta^{(2)}}^{-1}$. If $l=2$, then $\beta^{(3)}\sim \alpha^{(2)}$. However, $\alpha^{(2)}$ and $\beta^{(3)}$ belong to different tiles, a contradiction. If $l>2$, we can rewrite $\alpha^{(2)}=\alpha^{(21)}\circ \alpha^{(22)}$ and $\beta^{(2)}=\beta^{(21)}\circ \beta^{(22)}$, where $\alpha^{(21)}, \beta^{(21)}$ are of length $1$. Putting all of these together, we will obtain that $\alpha^{(22)}\vsim {\alpha^{(22)}}^{-1}$, which contradicts to Lemma \ref{l:nsim}. This completes the proof of $(2)$.
\end{proof}

The following useful result is also a direct consequence of the definition of equivalence and the compatibility of permissible arc segments.
\begin{lemma}
    Let $\eta$, $\gamma$ and $\delta$ be arc segments of length $l$ of certain pairwise compatible permissible arcs and denote by $\eta=\eta_1\circ\cdots\circ\eta_l$, $\gamma=\gamma_1\circ\cdots\circ \gamma_l$ and $\delta=\delta_1\circ\cdots\circ \delta_l$, where $\eta_i,\gamma_j,\delta_k$ are irreducible arc segments. Assume that for each $1\leq i\leq l$, $\eta_i,\gamma_i,\delta_i$ belong to the same tile $\Delta_i$.
    Let $\bf{a}$ be an arc in $\fT$ and fix an orientation of $\bf{a}$. Assume that $\eta(0)={\bf{a}}(t_1)$,  $\gamma(0)={\bf{a}}(t_2)$
and $\delta(0)={\bf{a}}(t_3)$ are in the interior of  $\bf{a}$ such that $t_1< t_2< t_3$. If $\eta\ \vec{\sim}\ \delta$, then $\eta\ \vec{\sim}\ \gamma$.
\end{lemma}

\subsection{A local to global criterion}\label{ss:local-to-global}

Let $(\fS,\fM,\fT)$ be a tiling and  $\Delta$ a tile.
For a finite multiset  $\cM$ of permissible arcs,  denote by 
  \begin{itemize}
      \item $\Arcseg_{\Delta}(\cM)$ the multiset of  irreducible arc segments of arcs in $\cM$ lying in $\Delta$;
      \item $\Arcseg^l(\cM)$ the multiset of  all the  arc segments of length $l$ for arcs in $\cM$;
      \item $\Arcseg^{\leq l}(\cM)$ the multiset of  all the  arc segments whose length is less than or equal to $l$ for arcs in $\cM$.
      
  \end{itemize}
  Notice that both $\Arcseg_{\Delta}(\cM)$ and $\Arcseg^{\leq l}(\cM)$ are finite multisets. Moreover, there is a positive integer $t$ such that $\cM$ is a subset of $\Arcseg^{\leq t}(\cM)$.
The following result provides a local criterion to identify two finite multisets of  permissible arcs, which plays a key role in the proof of Theorem \ref{t:main-result-1}.
\begin{lemma}\label{l:key lemma}
Let $\cM$, $\cN$ be two finite  multisets of pairwise compatible permissible arcs. The following are equivalent:
\begin{itemize}
\item[(1)] $\cM=\cN$;
\item[(2)] $\Arcseg^1(\cM)=\Arcseg^1(\cN)$.
\item[(3)] $\Arcseg^{\leq l}(\cM)=\Arcseg^{\leq l}(\cN)$ for each positive integer $l$.
\end{itemize}
\end{lemma}
\begin{proof}
It is evident that ``$(1)\Rightarrow (2)$". We can deduce that ``$(3)\Rightarrow (1)$" by observing that there exists a positive integer $t$ for which $\cM\subset \Arcseg^{\leq t}(\cM)$.
In order  to show ``$(2)\Rightarrow (3)$",  it suffices to prove that $\Arcseg^{\leq l}(\cM)=\Arcseg^{\leq l}(\cN)$ implies $\Arcseg^{l+1}(\cM)=\Arcseg^{l+1}(\cN)$. 
Consider $\km$, a maximal sub-multiset of $\Arcseg^{l+1}(\cM)$,  such that there exists a sub-multiset $\kn$ of $\Arcseg^{l+1}(\cN)$ with $\km=\kn$. In other words, if $\Arcseg^{l+1}(\cM)\backslash \km\neq \emptyset$, then for any $\alpha\in \Arcseg^{l+1}(\cM)\backslash\km$ and $\beta\in \Arcseg^{l+1}(\cN)\backslash\kn$, we have $\alpha\nsim \beta$. We need to prove that $\Arcseg^{l+1}(\cM)\backslash \km=\emptyset=\Arcseg^{l+1}(\cN)\backslash\kn$. If not, let us assume that $\Arcseg^{l+1}(\cM)\backslash\km\neq \emptyset$ and let $\alpha_1$ be an arc segment in $\Arcseg^{l+1}(\cM)\backslash\km$. We will demonstrate that there is an arc segment $\beta\in \Arcseg^{l+1}(\cN)\backslash\kn$ such that $\alpha_1\sim \beta$, which contradicts to the maximality of $\km$.

Let $\ba_1,\ba_2$ be the arcs in $\fT$ which divide $\alpha_1$ into three arc segments $\alpha_1^{(1)}$, $\alpha_1^{(2)}$ and $\alpha_1^{(3)}$ in order, \ie, $\alpha_1=\alpha_1^{(1)}\circ\alpha_1^{(2)}\circ\alpha_1^{(3)}$,  where  $\alpha_1^{(1)},\alpha_1^{(3)}\in{\Arcseg}^1(\cM)$ and $\alpha_1^{(2)}\in{\Arcseg}^{l-1}(\cM)$. 
Notice that $\alpha_1^{(2)}\circ \alpha_1^{(3)}\in{\Arcseg}^{l}(\cM)$ and ${\Arcseg}^{ l}(\cM)={\Arcseg}^{l}(\cN)$, there is an arc segment $\beta_1\in {\Arcseg}^{l+1}(\cN)$ divided by $\ba_1$ and $\ba_2$ into three arc segments  $\beta_1^{(1)}\in {\Arcseg}^{1}(\cN)$, $\beta_1^{(2)}\in{\Arcseg}^{l-1}(\cN)$ and $\beta_1^{(3)}\in {\Arcseg}^{1}(\cN)$ such that $\beta_1=\beta_1^{(1)}\circ \beta_1^{(2)}\circ\beta_1^{(3)}$,  $\beta_1^{(2)}\circ \beta_1^{(3)}\vsim \alpha_1^{(2)}\circ \alpha_1^{(3)}$ and $\alpha_1^{(1)}\nsim \beta_1^{(1)}$. If $\beta_1\in \kn$, there is an arc segment $\bar{\alpha}_1=\bar{\alpha}_1^{(1)}\circ\bar{\alpha}_1^{(2)}\circ\bar{\alpha}_1^{(3)}\in \km$ such that $\beta_1\vsim\bar{\alpha}_1$, where $\bar{\alpha}_1^{(1)}\in \Arcseg^1(\cM)$, $\bar{\alpha}_1^{(2)}\in \Arcseg^{l-1}(\cM)$ and $\bar{\alpha}_1^{(3)}\in\Arcseg^1(\cM)$ (cf. Figure \ref{f:Existence of beta_1 not in kn}). 
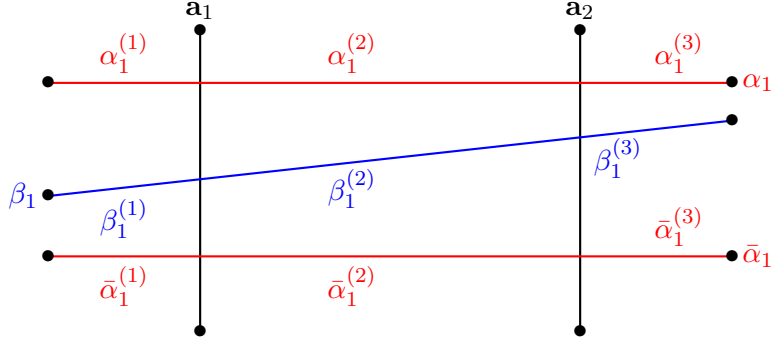
\begin{figure}[htpb]
		\begin{tikzpicture}
			

			\draw[black,thick] (-3,0)\nn to (-3,-4)\nn;
\draw[black,thick] (2,0)\nn to (2,-4)\nn;
 \draw[black,thick] (-3,0)node[above]{$\ba_1$}(2,0)node[above]{$\ba_2$}; 


     \draw[red,thick] (-5,-0.7)\nn to (4,-0.7)\nn;
						\draw[red] (4,-0.7)node[right]{$\alpha_1$};	
	\draw[red] (-4,-0.7)node[above]{$\alpha_1^{(1)}$} (-1,-0.7)node[above]{$\alpha_1^{(2)}$}(3.3,-0.7)node[above]{$\alpha_1^{(3)}$};

 \draw[blue,thick] (4,-1.2)\nn to (-5,-2.2)\nn;
			\draw[blue] (-5,-2.2)node[left]{$\beta_1$};
	\draw[blue] (-4,-2.1)node[below]{$\beta_1^{(1)}$} (-1,-2.5)node[above]{$\beta_1^{(2)}$}(2.5,-1.3)node[below]{$\beta_1^{(3)}$};

    \draw[red,thick] (-5,-3)\nn to (4,-3)\nn;
						\draw[red] (4,-3)node[right]{$\bar{\alpha}_1$};	
	\draw[red] (-4,-3)node[below]{$\bar{\alpha}_1^{(1)}$} (-1,-3)node[below]{$\bar{\alpha}_1^{(2)}$}(3.3,-3)node[above]{$\bar{\alpha}_1^{(3)}$}; 

		\end{tikzpicture}
		\caption{Existence of $\beta_1\not\in \kn$}\label{f:Existence of beta_1 not in kn}
	\end{figure}
We claim that $\alpha_1^{(2)}\circ\alpha_1^{(3)}\neq \bar{\alpha}_1^{(2)}\circ\bar{\alpha}_1^{(3)}$. Otherwise, we have  either  $(1):\ \alpha_1^{(2)}\circ \alpha_1^{(3)}\vsim \bar{\alpha}_1^{(2)}\circ \bar{\alpha}_1^{(3)}$ or  $(2):\ \alpha_1^{(2)}\circ \alpha_1^{(3)}\cev{\sim} \bar{\alpha}_1^{(2)}\circ \bar{\alpha}_1^{(3)}$. In case $(1)$,  we have $\alpha_1^{(1)}=\bar{\alpha}_1^{(1)}$. Hence $\alpha_1=\bar{\alpha}_1\in \km$, a contradiction. In case $(2)$, by applying Lemma \ref{l:equiv}, we obtain a contradiction.
Therefore $\alpha_1^{(2)}\circ\alpha_1^{(3)}\in \Arcseg^l(\cM)\backslash\{\bar{\alpha}_1^{(2)}\circ\bar{\alpha}_1^{(3)}\}=\Arcseg^l(\cN)\backslash\{\beta_1^{(2)}\circ\beta_1^{(3)}\}$. As a consequence, we may replace $\beta_1$ by a new arc segment. Since $\kn$ is finite, we may assume that $\beta_1\not\in \kn$.

Since $\beta_1^{(1)}\circ \beta_1^{(2)}\in \Arcseg^l(\cN)=\Arcseg^l(\cM)$. Similar to the previous discussion, there is an arc segment $\alpha_2\in \Arcseg^{l+1}(\cM)\backslash \km$ such that $\alpha_2=\alpha_2^{(1)}\circ\alpha_2^{(2)}\circ \alpha_2^{(3)}$, $\beta_1^{(1)}\circ \beta_1^{(2)}\vsim \alpha_2^{(1)}\circ \alpha_2^{(2)}$ and $\beta_1^{(3)}\nsim \alpha_2^{(3)}$ (cf. Figure~\ref{f:Existence of beta_1 and alpha_2}). 
\begin{figure}[htpb]
		\begin{tikzpicture}
			


			\draw[black,thick] (-3,0)\nn to (-3,-4)\nn;

\draw[black,thick] (2,0)\nn to (2,-4)\nn;

 \draw[black,thick] (-3,0)node[above]{$\ba_1$}(2,0)node[above]{$\ba_2$}; 


     \draw[red,thick] (-5,-0.7)\nn to (4,-0.7)\nn;
						\draw[red] (4,-0.7)node[right]{$\alpha_1$};	
	\draw[red] (-4,-0.7)node[above]{$\alpha_1^{(1)}$} (-1,-0.7)node[above]{$\alpha_1^{(2)}$}(3.3,-0.7)node[above]{$\alpha_1^{(3)}$} (-2.6,-0.7)node[above]{$P_1$};

 \draw[blue,thick] (4,-1.2)\nn to (-5,-2.2)\nn;
			\draw[blue] (-5,-2.2)node[left]{$\beta_1$};
	\draw[blue] (-4,-2.1)node[below]{$\beta_1^{(1)}$} (-1,-2.5)node[above]{$\beta_1^{(2)}$}(2.5,-1.3)node[below]{$\beta_1^{(3)}$}(-2.6,-1.9)node[below]{$Q_1$};

    \draw[red,thick] (-5,-3)\nn to (4,-3)\nn;
						\draw[red] (4,-3)node[right]{$\alpha_2$};	
	\draw[red] (-4,-3)node[below]{$\alpha_2^{(1)}$} (-1,-3)node[below]{$\alpha_2^{(2)}$}(3.3,-3)node[above]{$\alpha_2^{(3)}$};

		\end{tikzpicture}
		\caption{Existence of $\beta_1$ and $\alpha_2$}\label{f:Existence of beta_1 and alpha_2}
	\end{figure}
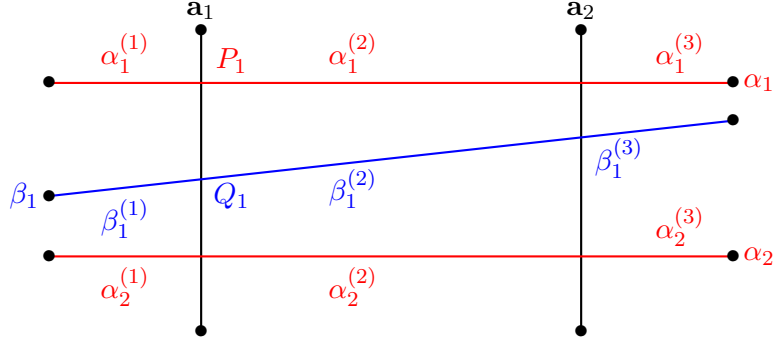

By repeating the procedure, we obtain arc segments
\[
\beta_1,\beta_2,\dots,\beta_t,\dots\in \Arcseg^{l+1}(\cN)\backslash\kn\ \text{and}\ 
\alpha_2,\alpha_3,\dots, \alpha_{t+1},\dots \in \Arcseg^{l+1}(\cM)\backslash\km.
\]
Since $\Arcseg^{l+1}(\cM)$ is finite, there will exist $1\leq s<t$ such that $\alpha_{s}=\alpha_{t}$. We may select a minimal $1\leq t$ that satisfies the condition where  there is an $1\leq s<t$ such that  $\alpha_{t}=\alpha_{s}$ and $\alpha_{s},\alpha_{s+1},\dots, \alpha_{t-1}$ are pairwise different. If necessary, we may  replace $\alpha_1$ with $\alpha_{s}$, allowing us to assume that $s=1$. Consequently, we have arc segments $\alpha_1,\alpha_2,\dots, \alpha_{t}\in \Arcseg^{l+1}(\cM)\backslash\km$ and $\beta_2,\dots, \beta_{t-1}\in \Arcseg^{l+1}(\cN)\backslash\kn$ (cf. Figure \ref{f:Existence of beta_{j+1} and alpha_{j+2}}). 
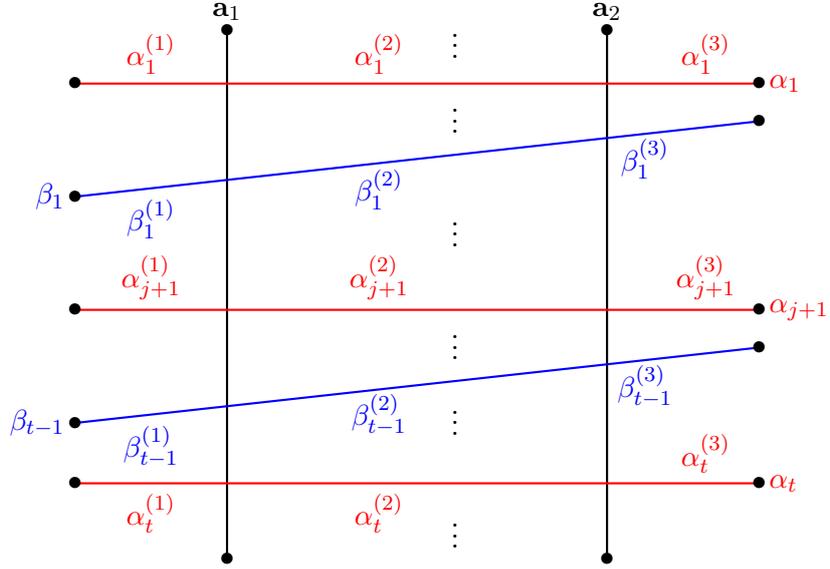
\begin{figure}[htpb]
		\begin{tikzpicture}
			



	\draw[black,thick] (-3,3)\nn to (-3,-4)\nn;

\draw[black,thick] (2,3)\nn to (2,-4)\nn;

 \draw[black,thick] (-3,3)node[above]{$\ba_1$}(2,3)node[above]{$\ba_2$}; 

	\draw[black,thick] (0,-1.5)node[above]{$\vdots$}; 
 \draw[black,thick] (0,2.5)node[above]{$\vdots$}; 
 \draw[black,thick] (0,-2.5)node[above]{$\vdots$}; 
 \draw[black,thick] (0,-4)node[above]{$\vdots$}; 
\draw[black,thick] (0,1.5)node[above]{$\vdots$}; 
\draw[black,thick] (0,0)node[above]{$\vdots$}; 
     \draw[red,thick] (-5,2.3)\nn to (4,2.3)\nn;
						\draw[red] (4,2.3)node[right]{$\alpha_1$};	
	\draw[red] (-4,2.3)node[above]{$\alpha_1^{(1)}$} (-1,2.3)node[above]{$\alpha_1^{(2)}$}(3.3,2.3)node[above]{$\alpha_1^{(3)}$}; 


 \draw[blue,thick] (4,1.8)\nn to (-5,0.8)\nn;
			\draw[blue] (-5,0.8)node[left]{$\beta_1$};
	\draw[blue] (-4,0.9)node[below]{$\beta_1^{(1)}$} (-1,0.5)node[above]{$\beta_1^{(2)}$}(2.5,1.7)node[below]{$\beta_1^{(3)}$};

     \draw[red,thick] (-5,-0.7)\nn to (4,-0.7)\nn;
						\draw[red] (4,-0.7)node[right]{$\alpha_{j+1}$};	
	\draw[red] (-4,-0.7)node[above]{$\alpha_{j+1}^{(1)}$} (-1,-0.7)node[above]{$\alpha_{j+1}^{(2)}$}(3.3,-0.7)node[above]{$\alpha_{j+1}^{(3)}$}; 


 \draw[blue,thick] (4,-1.2)\nn to (-5,-2.2)\nn;
			\draw[blue] (-5,-2.2)node[left]{$\beta_{t-1}$};
	\draw[blue] (-4,-2.1)node[below]{$\beta_{t-1}^{(1)}$} (-1,-2.5)node[above]{$\beta_{t-1}^{(2)}$}(2.5,-1.3)node[below]{$\beta_{t-1}^{(3)}$};

    \draw[red,thick] (-5,-3)\nn to (4,-3)\nn;
						\draw[red] (4,-3)node[right]{$\alpha_{t}$};	
	\draw[red] (-4,-3)node[below]{$\alpha_{t}^{(1)}$} (-1,-3)node[below]{$\alpha_{t}^{(2)}$}(3.3,-3)node[above]{$\alpha_{t}^{(3)}$};

		\end{tikzpicture}
		\caption{Existence of $\beta_j$ and $\alpha_{j+1}$}\label{f:Existence of beta_{j+1} and alpha_{j+2}}
	\end{figure}

For each $\alpha_i$ ($1\leq i\leq t$) and $\beta_j$ ($1\leq j<t$), we write
\[
\alpha_i=\alpha_i^{(1)}\circ \alpha_i^{(2)}\circ \alpha_i^{(3)}\ \text{and}\ \beta_j=\beta_j^{(1)}\circ \beta_j^{(2)}\circ \beta_j^{(3)}.
\]
By construction, we have
\begin{eqnarray*}
    &&\alpha_i^{(2)}\vsim \beta_j^{(2)}, 1\leq i\leq t, 1\leq j<t;\\
   && \alpha_i^{(2)}\circ \alpha_i^{(3)}\vsim \beta_{i}^{(2)}\circ \beta_i^{(3)}, 1\leq i<t;\\
    &&\beta_j^{(1)}\circ \beta_j^{(2)}\vsim \alpha_{j+1}^{(1)}\circ \alpha_{j+1}^{(2)}, 1\leq j<t;\\
    &&\alpha_i^{(1)}\nsim \beta_i^{(1)}, 1\leq i<t;\\
   && \beta_j^{(3)}\nsim\alpha_{j+1}^{(3)}, 1\leq j<t.
\end{eqnarray*}
Since $\alpha_1=\alpha_t$, we have $\alpha_1\vsim \alpha_t$  by Lemma \ref{l:equiv}. In particular, $\alpha_1^{(1)}=\alpha_t^{(1)}$. Consequently, $\alpha_1^{(1)}\vsim \beta_{t-1}^{(1)}$. More specifically, there exists an integer $t_0$ such that $1<t_0\leq t-1$ and $\beta_{t_0}^{(1)}\vsim \alpha_1^{(1)}$. Let $1<k\leq t-1$ be the smallest integer such that $\beta_k^{(1)}\vsim \alpha_1^{(1)}$.

We claim that $\alpha_1\vsim \beta_k$.
It is important to note that the arcs in $\cM$ (resp. $\cN$) are pairwise compatible. Hence arc segments in $\Arcseg^{\leq l}(\cM)$ (resp. $\Arcseg^{\leq l}(\cN)$) are pairwise compatible for each positive integer $l$. By the assumption that $\Arcseg^{\leq l}(\cM)=\Arcseg^{\leq l}(\cN)$, we conclude that every arc segment in $\Arcseg^{\leq l}(\cM)$  is compatible with the arc segments in $\Arcseg^{\leq l}(\cN)$.


 Fix an orientation of $\ba_1$. For two points $P=\ba_1(t_1)$ and $Q=\ba_1(t_2)$, we define $P<Q$ if $t_1<t_2$. 
The connection point of $\alpha_i^{(1)}$ and  $\alpha_i^{(2)}$ are denoted by $P_i$, and the connection point of $\beta_i^{(1)}$ and  $\beta_i^{(2)}$ are denoted by $Q_i$ for each $1\leq i\leq k$. 
Without loss of generality, we may suppose that $Q_1>P_1$.

Now, let us consider the case where $Q_k>P_1$. By $\beta_k^{(1)}\vsim\alpha_1^{(1)}$ and  $\beta_1^{(1)}\nsim \alpha_1^{(1)}$, we infer that $P_1<Q_k<Q_1$ (cf. Figure~\ref{f:case of qk>p1}). Since  $\beta_1^{(2)}\circ \beta_1^{(3)}\vsim \alpha_1^{(2)}\circ \alpha_1^{(3)}$ and $\beta_1^{(2)}\circ \beta_1^{(3)},\alpha_1^{(2)}\circ \alpha_1^{(3)}, \beta_k^{(2)}\circ \beta_k^{(3)}$ are compatible with each other, it follows that $ \beta_k^{(2)}\circ \beta_k^{(3)}\vsim \alpha_1^{(2)}\circ \alpha_1^{(3)}$. Consequently, we have $\beta_k=\beta_k^{(1)} \circ \beta_k^{(2)}\circ \beta_k^{(3)}\vsim \alpha_1^{(1)}\circ \alpha_1^{(2)}\circ \alpha_1^{(3)}=\alpha_1$.
\begin{figure}[htpb]
		\begin{tikzpicture}
			

			\draw[black,thick,->-=.5,>=stealth] (-3,0)\nn to (-3,-5)\nn;

\draw[black,thick] (2,0)\nn to (2,-5)\nn;

 \draw[black,thick] (-3,0)node[above]{$\ba_1$}(2,0)node[above]{$\ba_2$};

   \draw[black,thick] (0,-0.5)node[above]{$\vdots$}; 
	\draw[black,thick] (0,-1.5)node[above]{$\vdots$}; 
 \draw[black,thick] (0,-2.4)node[above]{$\vdots$}; 
 \draw[black,thick] (0,-3.5)node[above]{$\vdots$}; 
 \draw[black,thick] (0,-5)node[above]{$\vdots$}; 

     \draw[red,thick] (-5,-0.7)\nn to (4,-0.7)\nn;
						\draw[red] (4,-0.7)node[right]{$\alpha_1$};	
	\draw[red] (-4,-0.7)node[above]{$\alpha_1^{(1)}$} (-1,-0.7)node[above]{$\alpha_1^{(2)}$}(3.3,-0.7)node[above]{$\alpha_1^{(3)}$} (-2.6,-0.7)node[above]{$P_1$}; 

\draw[black,thick] (0,-0.5)node[above]{$\vdots$};

 \draw[blue,thick] (4,-1.5)\nn to (-5,-3.5)\nn;
			\draw[blue] (-5,-3.5)node[left]{$\beta_1$};
	\draw[blue] (-4,-3.2)node[below]{$\beta_1^{(1)}$} (-1,-3.4)node[above]{$\beta_1^{(2)}$}(3,-1.6)node[below]{$\beta_1^{(3)}$}(-2.6,-3)node[below]{$Q_1$}; 	

    \draw[red,thick] (-5,-4)\nn to (4,-4)\nn;
						\draw[red] (4,-4)node[right]{$\alpha_3$};	
	\draw[red] (-4,-4)node[below]{$\alpha_2^{(1)}$} (-1,-4)node[below]{$\alpha_2^{(2)}$}(3.3,-4)node[above]{$\alpha_2^{(3)}$};

 \draw[blue,thick] (4,-1.2)\nn to (-5,-2.2)\nn;
			\draw[blue] (-5,-2.2)node[left]{$\beta_k$};
	\draw[blue] (-4,-2.1)node[below]{$\beta_k^{(1)}$} (-1,-2.5)node[above]{$\beta_k^{(2)}$}(2.5,-0.6)node[below]{$\beta_k^{(3)}$}(-2.6,-1.9)node[below]{$Q_k$};

		\end{tikzpicture}
		\caption{Case of $Q_k>P_1$}\label{f:case of qk>p1}
	\end{figure}
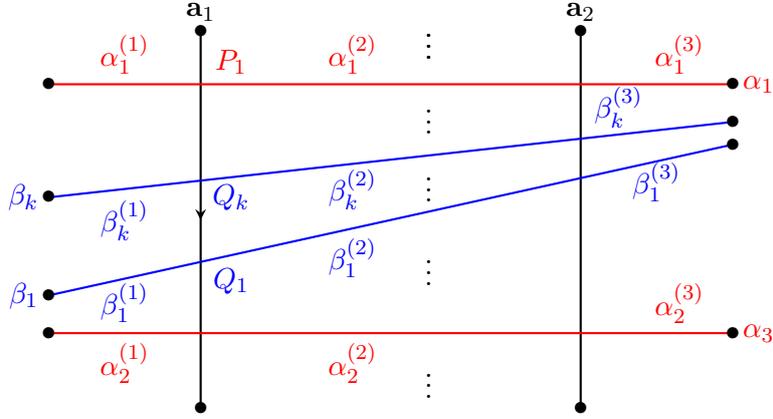

Now suppose that $Q_k<P_1$. 

If $P_k>P_1$,  we conclude that $\beta_k^{(2)}\circ \beta_k^{(3)}\vsim \alpha_1^{(2)}\circ \alpha_1^{(3)}$ by noticing that $\beta_k^{(2)}\circ \beta_k^{(3)}\vsim \alpha_k^{(2)}\circ \alpha_k^{(3)}$ and $\beta_k^{(2)}\circ \beta_k^{(3)},  \alpha_k^{(2)}\circ \alpha_k^{(3)}, \alpha_1^{(2)}\circ \alpha_1^{(3)}$ are compatible with each other.  As a consequence, $\beta_k=\beta_k^{(1)} \circ \beta_k^{(2)}\circ \beta_k^{(3)}\vsim \alpha_1^{(1)}\circ \alpha_1^{(2)}\circ \alpha_1^{(3)}=\alpha_1$ (cf. Figure~\ref{f:case of pk>p1>qk}).
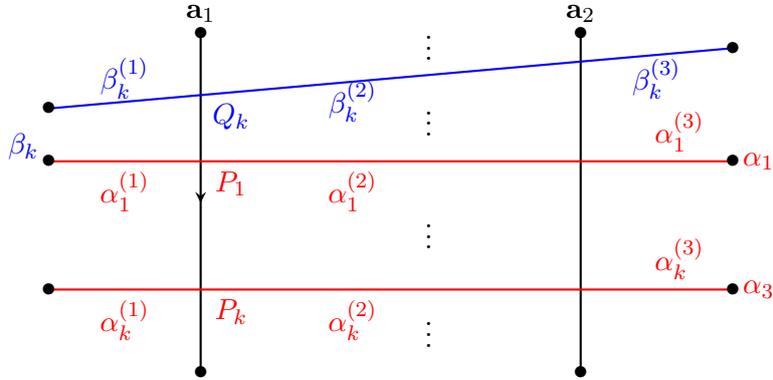
\begin{figure}[htpb]
		\begin{tikzpicture}
			

			\draw[black,thick,->-=.5,>=stealth] (-3,1)\nn to (-3,-3.5)\nn;


\draw[black,thick] (2,1)\nn to (2,-3.5)\nn;

 \draw[black,thick] (-3,1)node[above]{$\ba_1$}(2,1)node[above]{$\ba_2$};

    \draw[black,thick] (0,0.5)node[above]{$\vdots$}; 
        \draw[black,thick] (0,-0.5)node[above]{$\vdots$}; 
 \draw[black,thick] (0,-2)node[above]{$\vdots$}; 
 \draw[black,thick] (0,-3.3)node[above]{$\vdots$};

 \draw[blue,thick] (4,0.8)\nn to (-5,0)\nn;
				\draw[blue] (-5,-0.5)node[left]{$\beta_k$};
	\draw[blue] (-4,0)node[above]{$\beta_k^{(1)}$} (-1,0.5)node[below]{$\beta_k^{(2)}$}(3,0.8)node[below]{$\beta_k^{(3)}$}(-2.6,0.2)node[below]{$Q_k$};

     \draw[red,thick] (-5,-0.7)\nn to (4,-0.7)\nn;
						\draw[red] (4,-0.7)node[right]{$\alpha_1$};	
	\draw[red] (-4,-0.7)node[below]{$\alpha_1^{(1)}$} (-1,-0.7)node[below]{$\alpha_1^{(2)}$}(3.3,-0.7)node[above]{$\alpha_1^{(3)}$}(-2.6,-0.7)node[below]{$P_1$};

    \draw[red,thick] (-5,-2.4)\nn to (4,-2.4)\nn;
						\draw[red] (4,-2.4)node[right]{$\alpha_3$};	
	\draw[red] (-4,-2.4)node[below]{$\alpha_k^{(1)}$} (-1,-2.4)node[below]{$\alpha_k^{(2)}$}(3.3,-2.4)node[above]{$\alpha_k^{(3)}$}(-2.6,-2.4)node[below]{$P_k$};

		\end{tikzpicture}
		\caption{Case of $Q_k<P_1$ and $P_k>P_1$}\label{f:case of pk>p1>qk}
	\end{figure}

Now assume that $P_k<P_1$.
Let $j$ be the smallest integer such that $Q_j<P_1$. This implies that $2\leq j\leq k$ and  $Q_{j-1}>P_1$. We claim that $P_j>P_1$. Otherwise, $P_j<P_1$, then  $P_{j}<P_1<Q_{j-1}$. Since $\beta_{j-1}^{(1)}$, $\alpha_{j}^{(1)}$ and $\alpha_1^{(1)}$ are compatible and $\beta_{j-1}^{(1)}\sim\alpha_{j}^{(1)}$ (cf. Figure~\ref{f:case of P_j<P_1<Q_{j-1}}),  we deduce that $\beta_{j-1}^{(1)}\sim\alpha_{j}^{(1)}\sim\alpha_1^{(1)}$,  which contradicts to that $k$ is the smallest integer such that $\beta_k^{(1)}\sim\alpha_1^{(1)}$.
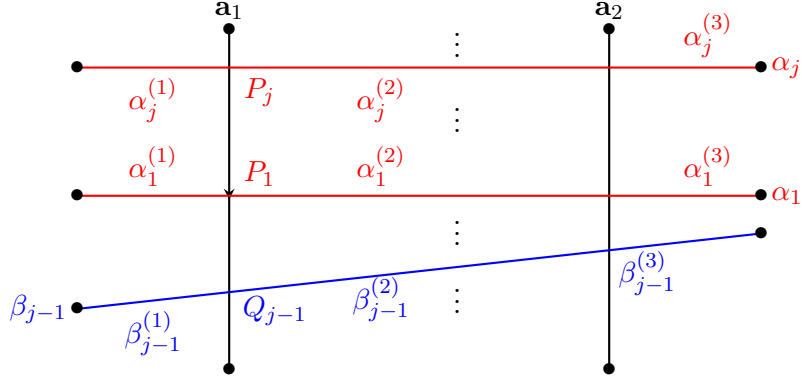
\begin{figure}[htpb]
		\begin{tikzpicture}
			

			\draw[black,thick,->-=.5,>=stealth] (-3,1.5)\nn to (-3,-3)\nn;


\draw[black,thick] (2,1.5)\nn to (2,-3)\nn;
 \draw[black,thick] (-3,1.5)node[above]{$\ba_1$}(2,1.5)node[above]{$\ba_2$};

   \draw[black,thick] (0,0)node[above]{$\vdots$}; 
	\draw[black,thick] (0,-1.5)node[above]{$\vdots$}; 
 \draw[black,thick] (0,-2.4)node[above]{$\vdots$}; 
 \draw[black,thick] (0,1)node[above]{$\vdots$}; 

    \draw[red,thick] (-5,1)\nn to (4,1)\nn;
						\draw[red] (4,1)node[right]{$\alpha_j$};	
	\draw[red] (-4,1)node[below]{$\alpha_j^{(1)}$} (-1,1)node[below]{$\alpha_j^{(2)}$}(3.3,1)node[above]{$\alpha_j^{(3)}$}(-2.6,1)node[below]{$P_j$};

     \draw[red,thick] (-5,-0.7)\nn to (4,-0.7)\nn;
						\draw[red] (4,-0.7)node[right]{$\alpha_1$};	
	\draw[red] (-4,-0.7)node[above]{$\alpha_1^{(1)}$} (-1,-0.7)node[above]{$\alpha_1^{(2)}$}(3.3,-0.7)node[above]{$\alpha_1^{(3)}$} (-2.6,-0.7)node[above]{$P_1$};

 \draw[blue,thick] (4,-1.2)\nn to (-5,-2.2)\nn;
			\draw[blue] (-5,-2.2)node[left]{$\beta_{j-1}$};
	\draw[blue] (-4,-2.1)node[below]{$\beta_{j-1}^{(1)}$} (-1,-2.5)node[above]{$\beta_{j-1}^{(2)}$}(2.5,-1.3)node[below]{$\beta_{j-1}^{(3)}$}(-2.4,-1.9)node[below]{$Q_{j-1}$};

		\end{tikzpicture}
		\caption{Case of $P_j<P_1<Q_{j-1}$}\label{f:case of P_j<P_1<Q_{j-1}}
	\end{figure}
So $P_j>P_1$. Because $\beta_j^{(2)}\circ \beta_j^{(3)}\vsim \alpha_j^{(2)}\circ \alpha_j^{(3)}$ and $\beta_j^{(2)}\circ \beta_j^{(3)},  \alpha_j^{(2)}\circ \alpha_j^{(3)}, \alpha_1^{(2)}\circ \alpha_1^{(3)}$ are compatible with each other, then $\beta_j^{(2)}\circ \beta_j^{(3)}\vsim \alpha_1^{(2)}\circ \alpha_1^{(3)}$. If $j=k$, then  $\beta_k=\beta_k^{(1)} \circ \beta_k^{(2)}\circ \beta_k^{(3)}\vsim \alpha_1^{(1)}\circ \alpha_1^{(2)}\circ \alpha_1^{(3)}=\alpha_1$. If $j<k$, then $\beta_j^{(1)}\nsim \alpha_1^{(1)}$. Because 
$\beta_k^{(1)}\sim \alpha_1^{(1)}$, we have $P_1>Q_k>Q_j$. Since $\beta_j^{(2)}\circ \beta_j^{(3)}\vsim \alpha_1^{(2)}\circ \alpha_1^{(3)}$ 
and $\beta_j^{(2)}\circ \beta_j^{(3)}, \alpha_1^{(2)}\circ \alpha_1^{(3)}, \beta_k^{(2)}\circ \beta_k^{(3)}$ are compatible with each other, we conclude that $\beta_k^{(2)}\circ \beta_k^{(3)}\vsim \alpha_1^{(2)}\circ \alpha_1^{(3)}$. As a result, $\beta_k=\beta_k^{(1)} \circ \beta_k^{(2)}\circ \beta_k^{(3)}\vsim \alpha_1^{(1)}\circ \alpha_1^{(2)}\circ \alpha_1^{(3)}=\alpha_1$ (cf. Figure~\ref{f:case of pk<p1>qk}). This completes the proof.
\begin{figure}[htpb]
		\begin{tikzpicture}
			

			\draw[black,thick,->-=.5,>=stealth] (-3,2)\nn to (-3,-3)\nn;
\draw[black,thick] (2,2)\nn to (2,-3)\nn;
 \draw[black,thick] (-3,2)node[above]{$\ba_1$}(2,2)node[above]{$\ba_2$};

 \draw[black,thick] (0,1.5)node[above]{$\vdots$}; 
    \draw[black,thick] (0,0.5)node[above]{$\vdots$}; 
     \draw[black,thick] (0,-0.5)node[above]{$\vdots$}; 
 \draw[black,thick] (0,-2)node[above]{$\vdots$}; 
 \draw[black,thick] (0,-3)node[above]{$\vdots$};

 \draw[blue,thick] (4,1.8)\nn to (-5,1)\nn;
				\draw[blue] (-5,1)node[left]{$\beta_j$};
	\draw[blue] (-4,1)node[above]{$\beta_j^{(1)}$} (-1,1.5)node[below]{$\beta_j^{(2)}$}(3,1.8)node[below]{$\beta_j^{(3)}$}(-2.6,1.2)node[below]{$Q_j$};

 \draw[blue,thick] (4,0.8)\nn to (-5,0)\nn;
				\draw[blue] (-5,-0.5)node[left]{$\beta_k$};
	\draw[blue] (-4,0)node[above]{$\beta_k^{(1)}$} (-1,0.5)node[below]{$\beta_k^{(2)}$}(3,0.8)node[below]{$\beta_k^{(3)}$}(-2.6,0.2)node[below]{$Q_k$};

     \draw[red,thick] (-5,-0.7)\nn to (4,-0.7)\nn;
						\draw[red] (4,-0.7)node[right]{$\alpha_1$};	
	\draw[red] (-4,-0.7)node[below]{$\alpha_1^{(1)}$} (-1,-0.7)node[below]{$\alpha_1^{(2)}$}(3.3,-0.7)node[above]{$\alpha_1^{(3)}$}(-2.6,-0.7)node[below]{$P_1$};

    \draw[red,thick] (-5,-2.4)\nn to (4,-2.4)\nn;
						\draw[red] (4,-2.4)node[right]{$\alpha_j$};	
	\draw[red] (-4,-2.4)node[below]{$\alpha_j^{(1)}$} (-1,-2.4)node[below]{$\alpha_j^{(2)}$}(3.3,-2.4)node[above]{$\alpha_j^{(3)}$}(-2.6,-2.4)node[below]{$P_j$};

		\end{tikzpicture}
		\caption{Case of $Q_k<P_1$, $P_k<P_1$ and $P_j>P_1$}\label{f:case of pk<p1>qk}
	\end{figure}
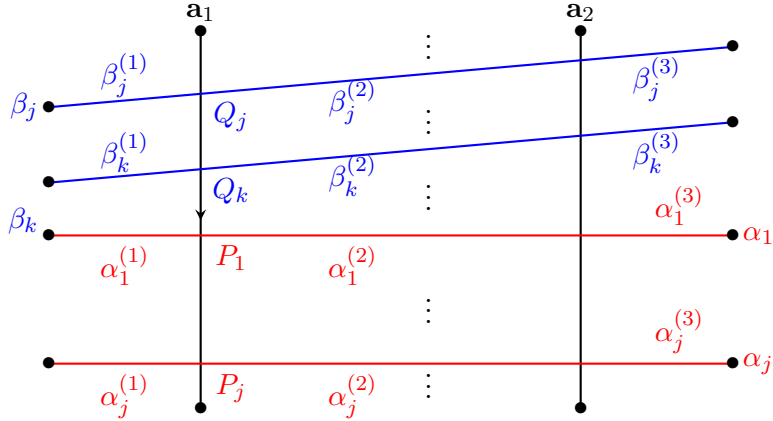

\end{proof}

\section{Intersection vectors}\label{s:intersection-vector}
\subsection{Intersection vectors}
Let $(\fS,\fM,\fT)$ be a tiling.  Recall that $\mathscr{R}(\fT)$ is the set of finite multisets of pairwise compatible permissible arcs.

Let $\beta$ be an arc on $\fS$ and  $\cM\in \mathscr{R}(\fT)$, set
\[\Int(\beta,\cM)=\Int(\cM,\beta):=\sum_{\gamma\in \cM}\Int(\gamma, \beta).
\]
The {\it intersection vector} $\Intv_\fT(\beta)$ of $\beta$ with respect to $\fT$ is defined as
\[\Intv_\fT(\beta)=(\Int(\beta,\ba_i))_{\ba_i\in \fT}.
\]
The {\it intersection vector} $\Intv_\fT(\cM)$ of $\cM$ with respect to $\fT$ is defined as
\[\Intv_\fT(\cM):=\sum_{\gamma\in \cM}\Intv_\fT(\gamma).
\]
It is clear that $\Intv_\fT(\cM)=(\Int(\cM,\ba_i))_{\ba_i\in \fT}$.

\subsection{Notations}
 Let $(\fS,\fM,\fT)$ be a tiling and $\Delta$ a tile. For each vertex $P$  of $\Delta$, there exists a unique edge  of $\Delta$, denoted by $E^{\Delta}_P$,  such that the arc segment connecting $P$ and one interior point of $E^{\Delta}_P$ is permissible.  For a finite multiset  $\cM\in \mathscr{R}(\fT)$,  denote by $\Arcseg(\cM, P, E^{\Delta}_P)$ the multiset of irreducible arc segments of arcs in $\cM$ lying in $\Delta$   with one endpoint   $P$ and the other in the interior of $E^{\Delta}_P$ and by 
 $\Seg(\cM, P, E^{\Delta}_P)$ the size  of $\Arcseg(\cM, P, E^{\Delta}_P)$. 
For an angle   $\angle \eta$  of $\Delta$, denoted by $\Arcseg(\cM, \angle \eta)$ the multiset of irreducible arc segments of arcs in  $\cM$  crossing $\angle \eta$ and by $\Seg(\cM, \angle \eta)$ the size of  $\Arcseg(\cM, \angle \eta)$.

\subsection{Proof of Theorem \ref{t:main-result-1}}
The aim of this section is to prove Theorem \ref{t:main-result-1}. The strategy of the proof is to show that for each tile $\Delta$ of $(\fS,\fM,\fT)$, the condition $\Intv_{\fT}(\cM)=\Intv_{\fT}(\cN)$ implies the local identities in Lemma \ref{l:key lemma}. We break down the lengthy proof into several lemmas by discussing the types of $\Delta$.
From now until the end of this section, we assume that $\cM,\cN\in \mathscr{R}(\fT)$ with $\Intv_\fT(\cM)=\Intv_\fT(\cN)$.  In particular,  $\Int(\cM,\ba)=\Int(\cN,\ba)$ for each $\ba\in \fT$.

Note that $\cM,\cN\in \mathscr{R}(\fT)$,   each arc segment in $\Arcseg(\cM)$ and $\Arcseg(\cN)$ is permissible. It suffices to prove that \[\Arcseg(\cM, P, E^{\Delta}_P)=\Arcseg(\cN, P, E^{\Delta}_P)\ \text{and}\ 
\Arcseg(\cM, \angle \eta)=\Arcseg(\cM, \angle \eta)\] for each vertex $P$ and each angle $\angle \eta$ for any tile $\Delta$. Equivalently, we need to prove  \[\Seg(\cM, P, E^{\Delta}_P)=\Seg(\cN, P, E^{\Delta}_P)\ \text{and}\  
\Seg(\cM, \angle \eta)=\Seg(\cM, \angle \eta).\]

\subsubsection{$\Delta$ is of type $(I)$}
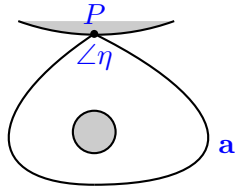
\begin{figure}[h]
\begin{tikzpicture}

\draw[thick,fill=black!20] (2,0.17)arc (250:290:3);

\fill(3,0) circle(1.5pt);

\draw[thick,color=black] (3,0) .. controls (1.5,-1) and (1.5,-2) .. (3,-2);
\draw[thick,color=black] (3,0) .. controls (5,-1) and (5,-2) .. (3,-2);

\draw[blue](3,0)node[above]{$P$}(4.5,-1.5)node[right]{$\ba$}(3,0)node[below]{$\angle \eta$};
\draw[thick,fill=black!20] (3,-1.3) circle(8pt);

\end{tikzpicture}
\caption{$\Delta$ is of Type I}\label{f:tile-type 1}
\end{figure}

Let  $\Delta$ be a tile of type $\text{(I)}$, $P$ the unique vertex,  and $\angle \eta$ the unique angle of $\Delta$ (cf. Figure~\ref{f:tile-type 1}). 

\begin{lemma}\label{l: one loop}
$\Seg(\cM, \angle \eta)=\Seg(\cN, \angle \eta)$  and  $\Seg(\cM, P, E^{\Delta}_P)= \Seg(\cN, P, E^{\Delta}_P)=0$.
\end{lemma}
\begin{proof}
Since $\Delta$ is of type $\text{(I)}$, $E^{\Delta}_P$ is the loop $\ba$.
Since each arc in $\cM$ is permissible and has no self-intersection, we have $\Seg(\cM, P, E^{\Delta}_P)=0$.  Similarly, $\Seg(\cN, P, E^{\Delta}_P)=0$.
Consequently,  $\Int(\cM,\ba)=2\Seg(\cM,\angle \eta)$ and $\Int(\cN,\ba)=2\Seg(\cN,\angle \eta)$. We conclude that $\Seg(\cM, \angle \eta)=\Seg(\cN, \angle \eta)$ by $\Int(\cM,\ba)=\Int(\cN,\ba)$. 
\end{proof}

\begin{figure}[h]

		\begin{tikzpicture}

			
		\draw[black,thick] (0,2)\nn to (-1,2)\nn;
			\draw[black,thick] (-1,2)\nn to (-2,1.5)\nn;
				\draw[black,thick] (-2,1.5)\nn to (-2.7,0.5)\nn;
			
				\draw[black,thick] (0,-2)\nn to (-1,-2)\nn;
			\draw[blue,thick] (-1,-2)\nn to (-2,-1.5)\nn;
				\draw[black,thick] (-2,-1.5)\nn to (-2.7,-0.5)\nn;
		
				\draw[black,thick] (0,2)\nn to (1,1.5)\nn;
			\draw[black,thick] (1,1.5)\nn to (1.7,0.5)\nn;
			
		\draw[black,thick] (0,-2)\nn to (1,-1.5)\nn;
			\draw[black,thick] (1,-1.5)\nn to (1.7,-0.5)\nn;

		\draw[blue](-2.7,-0.3)node[above]{$\vdots$}(1.7,-0.3)node[above]{$\vdots$}(0,2)node[above]{$P_1$}(-1,2)node[above]{$P_2$}(-2,1.5)node[above]{$P_3$}(-2.7,-0.5)node[left]{$P_{i-2}$}(1,1.5)node[above]{$P_m$}(1.7,0.5)node[right]{$P_{m-1}$}(-2,-1.5)node[left]{$P_{i-1}$}(-1,-2)node[below]{$P_{i}$};
		\draw[blue](-0.5,2)node[above]{$\gamma_1$}(-2.3,-1)node[left]{$\gamma_{i-2}$}(0.5,1.8)node[above]{$\gamma_m$}(-1.5,1.7)node[above]{$\gamma_2$}(-1.1,-2)node[left]{$\gamma_{i-1}$}(-0.5,-2)node[below]{$\gamma_{i}$};
		\draw[blue](-1.6,-1.5)node[above]{$\angle\eta_{i-1}$}(-1,-2)node[above]{$\angle\eta_{i}$}(0,2)node[below]{$\angle\eta_1$};


		\end{tikzpicture}

		\caption{ $ \Delta$ is an $m$-gons}\label{f:m-gons}
	\end{figure}
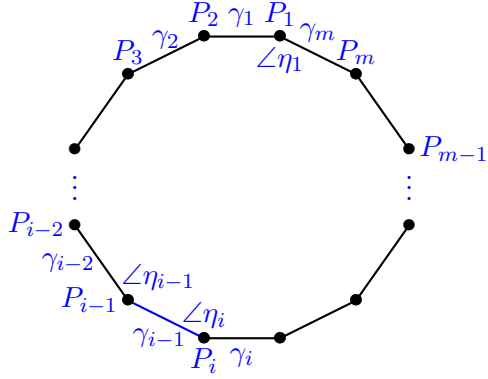
\subsubsection{$\Delta$ is an $m(\geq 3)$-gon}
Throughout this section, when we write notations with subscripts outside the domain $1,\dots,m$, we will implicitly assume identification modulo $m$.
Let  $\Delta$ be an $m$-gon  with arcs $\gamma_1,\gamma_2,\dots, \gamma_m$ in order and its interior contains no unmarked boundary segment. Denote by $P_i$ the incident endpoint of $\gamma_i$ and $\gamma_{i-1}$, and by $\angle\eta_i$ the angle in $\Delta$ with vertex $P_i$ (cf. Figure \ref{f:m-gons}). In particular, $E^{\Delta}_{P_{i-1}}=\gamma_{i}$ for each $i\in \mathbb{Z}$. Throughout this section, we fix the notations as in Figure \ref{f:m-gons}.

Since $\cM$ consists of compatible permissible arcs,  we have
\[\Seg(\cM, P_{i-1}, \gamma_{i})\circ \Seg(\cM,\angle \eta_{i-1})=0
\]
and
\[\Seg(\cM, P_i, \gamma_{i+1})\circ \Seg(\cM, P_{i-1}, \gamma_{i})=0
\]
for each $i\in \mathbb{Z}$.
Furthermore, if $\Seg(\cM, P_{i-1}, \gamma_{i})\neq 0$, then 
\[\Int(\gamma_{i-1},\cM)=\Seg(\cM,\angle \eta_{i})\] for each $i\in \mathbb{Z}$.

\begin{lemma}\label{l:huafenyizhi}
  For each $1\leq i\leq m$, $\Seg(\cM, P_{i-1}, \gamma_{i})\neq 0$ if and only if $\Seg(\cN, P_{i-1}, \gamma_{i})\neq 0$.
\end{lemma}
\begin{proof}
Let $S_1$ be the set of integers $i\in \{1,\dots,m\}$ such that $\Seg(\cM, P_{i-1}, \gamma_{i})\neq 0$ and  $S_2$ the set of integers $j\in \{1,\dots,m\}$ such that $\Seg(\cN, P_{j-1}, \gamma_{j})\neq 0$. It suffices to show that $S_1=S_2$. There is nothing to prove if both $S_1$ and $S_2$ are empty. Without loss of generality, we may assume that $S_1\neq \emptyset$ and denote by $S_1=\{i_1,i_2,\dots, i_c\}$, where $i_1<i_2<\cdots<i_c$. 
Denote by $i_{c+1}=i_1+m$ and set $\hat{S}_1=S_1\cup\{i_{c+1}\}$. For each $1\leq k\leq c$,
 we have $i_{k+1}-i_k\geq 2$.

We first show that $S_2\neq \emptyset$. Otherwise, assume that $S_2=\emptyset$. We define a sign function $\xi:\{\gamma_1,\dots, \gamma_m\}\to \{\pm 1\}$ as follows (cf. Figure~\ref{f:sign}):
\begin{itemize}
\item $\xi(\gamma_{i_k})=+1$ for $1\leq k\leq c$;
\item Let $1\leq k\leq c$. We define
\[\xi(\gamma_{i_k+a})=\begin{cases}
(-1)^a &\text{for $1\leq a\leq i_{k+1}-i_k-2$;}\\
-1& \text{for $a=i_{k+1}-i_k-1$.}
\end{cases} 
\]
\end{itemize} 
\begin{figure}[h]

		\begin{tikzpicture}

		\draw[black,thick] (0,0)\nn to (1,0)\nn to (2,0)\nn to (3,0)\nn to (4,0)\nn to (5,0)\nn;
		
		\draw[black,thick] (6,0)\nn to (7,0)\nn to (8,0)\nn;	
		
		\draw[blue](0,0)node[left]{$\cdots$}(5.1,0)node[right]{$\cdots$}(8,0)node[right]{$\cdots$};

  	\draw[blue](0.5,0)node[below]{$\gamma_{i_k-1}$}(1.5,0)node[below]{$\gamma_{i_k}$}(2.5,0)node[below]{$\gamma_{i_k+1}$};

	\draw[blue](6.5,0)node[below]{$\gamma_{i_{k+1}-1}$}(7.5,0)node[below]{$\gamma_{i_{k+1}}$};

	\draw[red](0.5,0)node[above]{$-$}(1.5,0)node[above]{$+$}(2.5,0)node[above]{$-$}(3.5,0)node[above]{$+$}(4.5,0)node[above]{$-$};

	\draw[red](6.5,0)node[above]{$-$}(7.5,0)node[above]{$+$};

		\end{tikzpicture}
\caption{ Sign function}\label{f:sign}
	\end{figure}
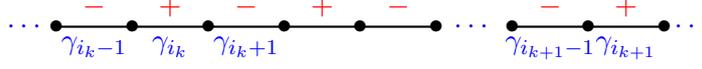
\begin{figure}[h]

		\begin{tikzpicture}

		\draw[black,thick] (0,0)\nn to (1,0)\nn to (2,0)\nn to (3,0)\nn to (4,0)\nn to (5,0)\nn;
		
		\draw[black,thick] (6,0)\nn to (7,0)\nn to (8,0)\nn;	
		
		\draw[blue](0,0)node[left]{$\cdots$}(5.1,0)node[right]{$\cdots$}(8,0)node[right]{$\cdots$};

  	\draw[blue](0.5,-0.2)node[below]{$\gamma_{i_k-1}$}(1.5,-0.2)node[below]{$\gamma_{i_k}$}(2.5,-0.2)node[below]{$\gamma_{i_k+1}$};

	\draw[blue](6.5,-0.2)node[below]{$\gamma_{i_{k+1}-1}$}(7.5,-0.2)node[below]{$\gamma_{i_{k+1}}$};

\draw[thick,color=red] (0,0) .. controls (1,1) and (1.5,1) .. (1.5,-0.2);
\draw[thick,color=red] (6,0) .. controls (7,1) and (7.5,1) .. (7.5,-0.2);

\draw[thick,color=red] (0.7,-0.2) .. controls (0.7,0.6) and (1.3,0.6) .. (1.3,-0.2);
\draw[thick,color=red] (1.7,-0.2) .. controls (1.7,0.6) and (2.3,0.6) .. (2.3,-0.2);
\draw[thick,color=red] (2.7,-0.2) .. controls (2.7,0.6) and (3.3,0.6) .. (3.3,-0.2);
\draw[thick,color=red] (3.7,-0.2) .. controls (3.7,0.6) and (4.3,0.6) .. (4.3,-0.2);
\draw[thick,color=red] (4.7,-0.2) .. controls (4.7,0.6) and (5.3,0.6) .. (5.3,-0.2);

\draw[thick,color=red] (6.7,-0.2) .. controls (6.7,0.6) and (7.3,0.6) .. (7.3,-0.2);

		\end{tikzpicture}
\caption{ Case for $\cM$}\label{f:case for m}
	\end{figure}
\begin{figure}[h]

		\begin{tikzpicture}

		\draw[black,thick] (0,0)\nn to (1,0)\nn to (2,0)\nn to (3,0)\nn to (4,0)\nn to (5,0)\nn;
		
		\draw[black,thick] (6,0)\nn to (7,0)\nn to (8,0)\nn ;	
		
		\draw[blue](0,0)node[left]{$\cdots$}(5.1,0)node[right]{$\cdots$}(8,0)node[right]{$\cdots$};

  	\draw[blue](0.5,-0.2)node[below]{$\gamma_{i_k-1}$}(1.5,-0.2)node[below]{$\gamma_{i_k}$}(2.5,-0.2)node[below]{$\gamma_{i_k+1}$};

	\draw[blue](6.5,-0.2)node[below]{$\gamma_{i_{k+1}-1}$}(7.5,-0.2)node[below]{$\gamma_{i_{k+1}}$};

\draw[thick,color=red] (-0.3,-0.2) .. controls (-0.3,0.6) and (0.3,0.6) .. (0.3,-0.2);

\draw[thick,color=red] (0.7,-0.2) .. controls (0.7,0.6) and (1.3,0.6) .. (1.3,-0.2);
\draw[thick,color=red] (1.7,-0.2) .. controls (1.7,0.6) and (2.3,0.6) .. (2.3,-0.2);
\draw[thick,color=red] (2.7,-0.2) .. controls (2.7,0.6) and (3.3,0.6) .. (3.3,-0.2);
\draw[thick,color=red] (3.7,-0.2) .. controls (3.7,0.6) and (4.3,0.6) .. (4.3,-0.2);
\draw[thick,color=red] (4.7,-0.2) .. controls (4.7,0.6) and (5.3,0.6) .. (5.3,-0.2);
\draw[thick,color=red] (5.7,-0.2) .. controls (5.7,0.6) and (6.3,0.6) .. (6.3,-0.2);
\draw[thick,color=red] (6.7,-0.2) .. controls (6.7,0.6) and (7.3,0.6) .. (7.3,-0.2);

		\end{tikzpicture}
\caption{ Case for $\cN$}\label{f:case for n}
	\end{figure}
 It follows that
\[\sum_{i=1}^m\xi(\gamma_i)\Int(\gamma_i,\cM)=\sum_{k=1}^c\xi(\gamma_{i_k})\Seg(\cM, P_{i_k-1},\gamma_{i_k})>0
\]
(cf. Figure~\ref{f:case for m}) and 
\[\sum_{i=1}^m\xi(\gamma_i)\Int(\gamma_i,\cN)\leq 0
\]
(cf. Figure~\ref{f:case for n}), which contradicts to \[\sum_{i=1}^m\xi(\gamma_i)\Int(\gamma_i,\cM)=\sum_{i=1}^m\xi(\gamma_i)\Int(\gamma_i,\cN).\]

Denote by $S_2=\{j_1,\dots, j_d\}$, where $j_1<j_2<\cdots<j_d$.
Suppose that $S_1\neq S_2$.
We separate the remaining proof by discussing  $S_1\cap S_2\neq \emptyset$ or $S_1\cap S_2=\emptyset$.

\noindent{\bf Case 1: $S_1\cap S_2\neq \emptyset$.} Denote by $j_{d+1}=j_1+m$ and $\hat{S}_2=S_2\cup\{j_{d+1}\}$. By $S_1\neq S_2$ and  $S_1\cap S_2\neq \emptyset$, we conclude that at least one of the following holds:
\begin{itemize}
\item[(1)] there is a consecutive pair $(i_k, i_{k+1})$ of $\hat{S}_1$ such that $\Seg(\cN, P_{i_{k+1}-1},\gamma_{i_{k+1}})\neq 0$ and $\Seg(\cN, P_{l-1}, \gamma_l)=0$ for $i_k\leq l<i_{k+1}$;
\item[(2)] there is a consecutive pair $(j_t, j_{t+1})$ of $\hat{S}_2$ such that  $\Seg(\cM, P_{j_{t+1}-1},\gamma_{j_{t+1}})\neq 0$ and $\Seg(\cM, P_{l-1}, \gamma_l)=0$ for $j_t\leq l<j_{t+1}$.
\end{itemize}
Without loss of generality, we assume that the condition $(1)$ holds and denote by $(a,b):=(i_k,i_{k+1})$ (cf. Figure \ref{f:ab}).

\begin{figure}[h]

		\begin{tikzpicture}

			
		\draw[black,thick] (0,2)\nn to (-1,2)\nn;
			\draw[black,thick] (-1,2)\nn to (-2,1.5)\nn;
				\draw[black,thick] (-2,1.5)\nn to (-2.7,0.5)\nn;
			
				\draw[black,thick] (0,-2)\nn to (-1,-2)\nn;
			\draw[black,thick] (-1,-2)\nn to (-2,-1.5)\nn;
				\draw[black,thick] (-2,-1.5)\nn to (-2.7,-0.5)\nn;
		
				\draw[black,thick] (0,2)\nn to (1,1.5)\nn;
			\draw[black,thick] (1,1.5)\nn to (1.7,0.5)\nn;
			
		\draw[black,thick] (0,-2)\nn to (1,-1.5)\nn;
			\draw[black,thick] (1,-1.5)\nn to (1.7,-0.5)\nn;

		\draw[blue](-2.7,-0.3)node[above]{$\vdots$}(1.7,-0.3)node[above]{$\vdots$}(0,2)node[above]{$P_b$}(-2.7,-0.5)node[left]{$P_{a-1}$}(1.7,0.5)node[right]{$P_{b-2}$}(-2,-1.5)node[below]{$P_{a}$}(-0.6,-2)node[below]{$P_{a+1}$};
		\draw[blue](-0.5,2)node[above]{$\gamma_b$}(-2.3,-1)node[left]{$\gamma_{a-1}$}(0.7,1.6)node[above]{$\gamma_{b-1}$}(-1,-1.6)node[left]{$\gamma_{a}$}(0,-1.6)node[left]{$\gamma_{a+1}$}(1.7,0.8)node[above]{$\gamma_{b-2}$};
		\draw[blue](-1.8,-1.5)node[above]{$\angle\eta_a$}(0,2)node[below]{$\angle\eta_b$};

\draw[thick,color=red] (-2.7,-0.5) .. controls (-1.5,-0.5) and (-1.5,-1) .. (-1.5,-2.4);
\draw[thick,color=red] (1,1.5) .. controls (0,1) and (-0.8,2) .. (-0.8,2.3);

\draw[thick,color=red] (-2.4,-1.7) .. controls (-2,-1) and (-1.7,-1.4) .. (-1.8,-2);
\draw[thick,color=red] (-1.4,-2.2) .. controls (-1.4,-1.5) and (-0.6,-1.5) .. (-0.6,-2.2);
\draw[thick,color=red] (-0.4,-2.2) .. controls (-0.4,-1.5) and (0.4,-1.5) .. (0.4,-2.2);
\draw[thick,color=red] (1.6,-1.7) .. controls (1,-1) and (0.7,-1.4) .. (0.8,-2);

\draw[thick,color=red] (-0.4,2.2) .. controls (-0.4,1.7) and (0.4,1.7) .. (0.4,2.2);

\draw[thick,color=red] (1.9,0.3) .. controls (1.4,0.3) and (1.4,0.7) .. (1.9,0.7);

\draw[thick,color=red] (1.9,-0.3) .. controls (1.4,-0.3) and (1.4,-0.7) .. (1.9,-0.7);

		\end{tikzpicture}
	
		\caption{ Case: $(a,b)=(i_k,i_{k+1})$ for $\cM$}\label{f:ab}
	\end{figure}

By definition of $S_1$, we have
\[\Seg(\cM, P_{b-1}, \gamma_{b})\neq 0\neq \Seg(\cM, P_{a-1}, \gamma_{a})
\]
and $\Seg(\cM, P_{l-1}, \gamma_{l})=0$  for each $a< l<b$.
As a result,
 \[
 \Int(\gamma_i,\cM)=
 \begin{cases}
 \Seg(\cM,\angle\eta_a)& i=a-1;\\
 \Seg(\cM,\angle\eta_a)+\Seg(\cM,P_{a-1},\gamma_a) &\text{$i=a$ and if $b-a=2$};\\
 \Seg(\cM,\angle\eta_{a})+\Seg(\cM,P_{a-1},\gamma_{a}) +\Seg(\cM,\angle\eta_{a+1})&\text{$i=a$ and if $b-a>2$};\\
 \Seg(\cM, \angle\eta_i)+\Seg(\cM, \angle\eta_{i+1}) & a+1\leq i\leq b-3;\\
 \Seg(\cM,\angle\eta_{b-2})& \text{$i=b-2$ and $b-a>2$.}
 \end{cases}
 \]
 
 In either case, we have
\[\sum_{i={a-1}}^{b-2}(-1)^i\Int(\gamma_{i},\cM)=(-1)^{a}\Seg(\cM,P_{a-1},\gamma_{a})\]
(cf. Figure~\ref{f:ab}). In particular, if $a$ is even, then $\sum_{i={a-1}}^{b-2}(-1)^i\Int(\gamma_{i},\cM)>0$, while if $a$ is odd, then $\sum_{i={a-1}}^{b-2}(-1)^i\Int(\gamma_{i},\cM)<0$.

On the other hand, we have
\[\Int(\gamma_{a-1},\cN)\ge \Seg(\cN,\angle\eta_{a}),\  \Int(\gamma_{b-2},\cN)=\Seg(\cN,\angle\eta_{b-2})\] and  for $a\leq i\leq b-3$,  \[\Int(\gamma_i,\cN)=\Seg(\cN,\angle\eta_i) +\Seg(\cN,\angle\eta_{i+1})\] (cf. Figure~\ref{f:ab for n}).  Hence
 \[\sum_{i={a-1}}^{b-2}(-1)^i\Int(\gamma_{i},\cN)=(-1)^{a-1}(\Int(\gamma_{a-1},\cN)-\Seg(\cN,\angle\eta_{a})).\]
 \begin{figure}[h]

		\begin{tikzpicture}

			
		\draw[black,thick] (0,2)\nn to (-1,2)\nn;
			\draw[black,thick] (-1,2)\nn to (-2,1.5)\nn;
				\draw[black,thick] (-2,1.5)\nn to (-2.7,0.5)\nn;
			
				\draw[black,thick] (0,-2)\nn to (-1,-2)\nn;
			\draw[black,thick] (-1,-2)\nn to (-2,-1.5)\nn;
				\draw[black,thick] (-2,-1.5)\nn to (-2.7,-0.5)\nn;
		
				\draw[black,thick] (0,2)\nn to (1,1.5)\nn;
			\draw[black,thick] (1,1.5)\nn to (1.7,0.5)\nn;
			
		\draw[black,thick] (0,-2)\nn to (1,-1.5)\nn;
			\draw[black,thick] (1,-1.5)\nn to (1.7,-0.5)\nn;

		\draw[blue](-2.7,-0.3)node[above]{$\vdots$}(1.7,-0.3)node[above]{$\vdots$}(0,2)node[above]{$P_b$}(-2.7,-0.5)node[left]{$P_{a-1}$}(1.7,0.5)node[right]{$P_{b-2}$}(-2,-1.5)node[below]{$P_{a}$}(-0.6,-2)node[below]{$P_{a+1}$};
		\draw[blue](-0.5,2)node[above]{$\gamma_b$}(-2.3,-1)node[left]{$\gamma_{a-1}$}(0.7,1.6)node[above]{$\gamma_{b-1}$}(-1,-1.6)node[left]{$\gamma_{a}$}(0,-1.6)node[left]{$\gamma_{a+1}$}(1.7,0.8)node[above]{$\gamma_{b-2}$};
		\draw[blue](-1.8,-1.5)node[above]{$\angle\eta_a$}(0,2)node[below]{$\angle\eta_b$}(1.6,0.5)node[left]{$\angle\eta_{b-2}$};

\draw[thick,color=red] (1,1.5) .. controls (0,1) and (-0.8,2) .. (-0.8,2.3);

\draw[thick,color=red] (-2.4,-1.7) .. controls (-2,-1) and (-1.7,-1.4) .. (-1.8,-2);
\draw[thick,color=red] (-1.4,-2.2) .. controls (-1.4,-1.5) and (-0.6,-1.5) .. (-0.6,-2.2);
\draw[thick,color=red] (-0.4,-2.2) .. controls (-0.4,-1.5) and (0.4,-1.5) .. (0.4,-2.2);
\draw[thick,color=red] (1.6,-1.7) .. controls (1,-1) and (0.7,-1.4) .. (0.8,-2);

\draw[thick,color=red] (-0.4,2.2) .. controls (-0.4,1.7) and (0.4,1.7) .. (0.4,2.2);

\draw[thick,color=red] (1.9,0.3) .. controls (1.4,0.3) and (1.4,0.7) .. (1.9,0.7);

\draw[thick,color=red] (1.9,-0.3) .. controls (1.4,-0.3) and (1.4,-0.7) .. (1.9,-0.7);

		\end{tikzpicture}
	
		\caption{ Case: $(a,b)=(i_k,i_{k+1})$ for  $\cN$}\label{f:ab for n}
	\end{figure}
Consequently, if $a$ is  even, then $\sum_{i={a-1}}^{b-2}(-1)^i\Int(\gamma_{i},\cN)\leq 0$, while if $a$ is odd, then $\sum_{i={a-1}}^{b-2}(-1)^i\Int(\gamma_{i},\cN)\ge 0$, which contradicts to \[\sum_{i={a-1}}^{b-2}(-1)^i\Int(\gamma_{i},\cM)
=\sum_{i={a-1}}^{b-2}(-1)^i\Int(\gamma_{i},\cN).\]

\noindent{\bf Case 2: $S_1\cap S_2=\emptyset$.}
 Define $j_{d+k}:=j_{k}+m$ for each $1\leq k\leq d$ and  $\tilde{S}_2:=\{j_1,j_2,\cdots,j_{2d}\}$.

We define a sign function $\xi:\{\gamma_1,\dots,\gamma_m\}\to \{\pm 1\}$ as follows:
\begin{itemize}
    \item $\xi(\gamma_{i_k})=+1$ for each $1\leq k\leq c$ and $\xi(\gamma_{j_k})=-1$ for $1\leq k\leq d$;
    \item for a consecutive pair $(i_k,i_{k+1})$ of $\hat{S}_1$, denote by $\mathcal{I}_k=\{j\in \tilde{S}_2~|~i_k<j<i_{k+1}\}$.
    \begin{itemize}
    \item If $\mathcal{I}_k=\emptyset$,  then define  (cf. Figure~\ref{f:sign})
    \[\xi(\gamma_{i_k+a})=\begin{cases}
        (-1)^a &\text{if $1\leq a\leq i_{k+1}-i_k-2$};\\
        -1& \text{if $a=i_{k+1}-i_k-1$}.
    \end{cases}
    \]
   
    \item Assume $\mathcal{I}_k\neq \emptyset$ and denote by $\mathcal{I}_k=\{j_{s_1},\dots, j_{s_t}\}$, where $j_{s_1}<j_{s_2}<\cdots<j_{s_t}$. We define (cf. Figure \ref{f:sign 2})
    \[\xi(\gamma_u)=\begin{cases}(-1)^a& \text{if $u=i_k+a$ and $1\leq a\leq j_{s_1}-i_k-2$};\\
    +1& \text{if $u=j_{s_1}-1$}.
    \end{cases}
    \]
    For  any $1\leq l<t$, we define
    \[\xi(\gamma_u)=\begin{cases}(-1)^{a+1}& \text{if $u=j_{s_l}+a$ and $1\leq a\leq j_{s_{l+1}}-j_{s_l}-2$};\\
    +1&\text{if $u=j_{s_{l+1}}-1$}.
    \end{cases}
    \]
    Finally, define
    \[\xi(\gamma_u)=\begin{cases}(-1)^{a+1}&\text{if $u=j_{s_t}+a$ and $1\leq a\leq i_{k+1}-j_{s_t}-2$};\\
    -1&\text{if $u=i_{k+1}-1$}.
    \end{cases}
    \]
    \end{itemize}
\end{itemize}

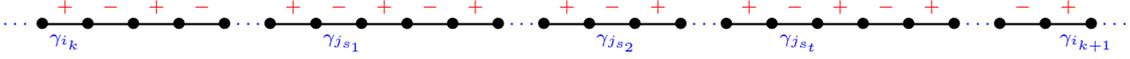
\begin{figure}[h]

		\begin{tikzpicture}[xscale=0.6,yscale=0.9]

		\draw[black,thick]  (1,0)\nn to (2,0)\nn to (3,0)\nn to (4,0) \nn to (5,0)\nn;
		
		\draw[black,thick]    (6,0)\nn to (7,0)\nn to (8,0)\nn to (9,0)\nn to  (10,0)\nn to  (11,0)\nn ;	

  \draw[black,thick]   (12,0)\nn to (13,0)\nn to (14,0)\nn to (15,0)\nn ;	

 \draw[black,thick](16,0)\nn to  (17,0)\nn to (18,0)\nn to (19,0)\nn to (20,0)\nn to (21,0)\nn ;

   \draw[black,thick]  (22,0)\nn to (23,0)\nn to (24,0)\nn;
		\draw[blue](1,0)node[left]{\tiny$\cdots$}(5,0)node[right]{\tiny$\cdots$}(11,0)node[right]{\tiny$\cdots$}(15,0)node[right]{\tiny$\cdots$}(21,0)node[right]{\tiny$\cdots$}
  (24,0)node[right]{\tiny$\cdots$};

  	\draw[blue](1.5,0)node[below]{\tiny$\gamma_{i_k}$};

	\draw[blue](7.6,0)node[below]{\tiny $\gamma_{j_{s_1}}$}(13.6,0)node[below]{\tiny $\gamma_{j_{s_2}}$}(17.6,0)node[below]{\tiny $\gamma_{j_{s_t}}$}(23.9,0)node[below]{\tiny $\gamma_{i_{k+1}}$};

	\draw[red](1.5,0)node[above]{\tiny$+$}(2.5,0)node[above]{\tiny$-$}(3.5,0)node[above]{\tiny$+$}(4.5,0)node[above]{\tiny$-$};

	\draw[red](6.5,0)node[above]{\tiny$+$}(7.5,0)node[above]{\tiny$-$} (8.5,0)node[above]{\tiny$+$}(9.5,0)node[above]{\tiny$-$}(10.5,0)node[above]{\tiny$+$}(12.5,0)node[above]{\tiny$+$}(13.5,0)node[above]{\tiny$-$}(14.5,0)node[above]{\tiny$+$}(16.5,0)node[above]{\tiny$+$}(17.5,0)node[above]{\tiny$-$}(18.5,0)node[above]{\tiny$+$}(20.5,0)node[above]{\tiny$+$}(19.5,0)node[above]{\tiny$-$}(22.5,0)node[above]{\tiny$-$}(23.5,0)node[above]{\tiny$+$};

		\end{tikzpicture}
\caption{ Sign function for case $\mathcal{I}_k\neq \emptyset$}\label{f:sign 2}
	\end{figure}
 It is easy to check that this is well-defined. For any $i_k<l\leq i_{k+1}$ such that $\Seg(\cM,\angle \eta_l)\neq 0$, we have $\xi(\gamma_{l-1})\xi(\gamma_l)=-1$ or $\xi(\gamma_{l-1})=\xi(\gamma_l)=+1$. It follows that
\[\sum_{i=1}^m\xi(\gamma_i)\Int(\gamma_i,\cM)\geq\sum_{i_k\in S_1}\xi(\gamma_{i_k})\Seg(\cM,P_{i_k-1},\gamma_{i_k})>0.\]
Similarly, for each $j_k<l\leq j_{k+1}$ such that $\Seg(\cN,\angle\eta_l)\neq 0$, we have $\xi(\gamma_{l-1})\xi(\gamma_l)=-1$ or $\xi(\gamma_{l-1})=\xi(\gamma_l)=-1$. Consequently, 
\[\sum_{i=1}^m\xi(\gamma_i)\Int(\gamma_i,\cN)\leq\sum_{j_k\in S_2}\xi(\gamma_{j_k})\Seg(\cN,P_{j_k-1},\gamma_{j_k})< 0,\]
which contradicts to 
\[\sum_{i=1}^m\xi(\gamma_i)\Int(\gamma_i,\cM)=\sum_{i=1}^m\xi(\gamma_i)\Int(\gamma_i,\cN).
\]
 This completes the proof.
\end{proof}

\begin{lemma}\label{l:m-gons}
Suppose that the tile $\Delta$ satisfies at least one of the following conditions:
\begin{itemize}
\item [$(1)$] $m$ is odd.
\item [$(2)$] $\Delta$ is of  type $(IV)$, \ie, $\Delta$ has exactly one boundary segment.
\end{itemize}
 For each $1\leq i\leq m$, we have \[\Seg(\cM, \angle \eta_i)=\Seg(\cN, \angle \eta_i), \ \ \ \ \Seg(\cM, P_i, \gamma_{i+1})= \Seg(\cN, P_i, \gamma_{i+1}).\]
\end{lemma}

\begin{proof}
Let $S$ be the subset of $\{1,\dots,m\}$ consisting of index $i$ such that $\Seg(\cM, P_{i-1},\gamma_i)\neq 0$. By Lemma \ref{l:huafenyizhi}, $\Seg(\cN, P_{i-1}, \gamma_i)\neq 0$ if and only if $i\in S$. Denote by $d_i:=\Int(\gamma_i,\cM)=\Int(\gamma_i,\cN)$ for $1\leq i\leq m$ for simplicity.

Let us first consider the case $S=\emptyset$. We have to show that $\Seg(\cM, \angle \eta_i)=\Seg(\cN, \angle \eta_i)$ for each $1\leq i\leq m$.
 In this case,  we have 
\[d_i=\Seg(\cM, \angle \eta_i)+\Seg(\cM, \angle\eta_{i+1})
\]
and 
\[d_i=\Seg(\cN, \angle\eta_i)+\Seg(\cN, \angle\eta_{i+1})
\]
for $1\leq i\leq m$.  We conclude that  $\Seg(\cM, \angle\eta_i)=\Seg(\cN,\angle\eta_i)$ for each $1\leq i\leq m$ by the fact that the equation \[
\begin{cases}x_1+x_2=d_1\\
\ \ \ \vdots\\
x_{m-1}+x_m=d_{m-1}\\
x_m+x_1=d_m\end{cases}
\] has a unique solution provided that $m$ is odd.

Assume that $\Delta$ is of type $\text{(IV)}$. Without loss of generality, we assume that $\gamma_1$ is the boundary segment. Hence we have
\[
\begin{cases}
\Seg(\cM,\angle\eta_1)=0\\
\Seg(\cM,\angle\eta_2)=0\\
\Seg(\cM,\angle\eta_3)=d_2\\
\Seg(\cM,\angle\eta_3)+\Seg(\cM,\angle\eta_4)=d_3\\
\ \ \ \vdots\\
\Seg(\cM,\angle\eta_{m-1})+\Seg(\cM,\angle\eta_m)=d_{m-1}\\
\Seg(\cM,\angle\eta_m)=d_m
\end{cases},\  \begin{cases}
\Seg(\cN,\angle\eta_1)=0\\
\Seg(\cN,\angle\eta_2)=0\\
\Seg(\cN,\angle\eta_3)=d_2\\
\Seg(\cN,\angle\eta_3)+\Seg(\cM,\angle\eta_4)=d_3\\
\ \ \ \vdots\\
\Seg(\cN,\angle\eta_{m-1})+\Seg(\cM,\angle\eta_m)=d_{m-1}\\
\Seg(\cN,\angle\eta_m)=d_m
\end{cases}.
\]
 As a consequence,  $\Seg(\cM,\angle\eta_i)=\Seg(\cN,\angle\eta_i)$ for $1\leq i\leq m$.
 
 Now assume that $S\neq \emptyset$ and denote by $S=\{i_1,\dots, i_c\}$, where $i_1<\cdots<i_c$. Set $i_{c+1}=i_c+m$ and denote by $\hat{S}=S\cup \{i_{c+1}\}$. If $\Delta$ is of type $\text{(IV)}$ with the boundary segment $\gamma_1$, we can view $\Delta$ as an $m$-gon without boundary segment by setting \[\Seg(\cM, \angle\eta_1)=0=\Seg(\cM,\angle\eta_2),\ \Seg(\cM, P_{m}, \gamma_1)=0\]
 \[\Seg(\cN, \angle\eta_1)=0=\Seg(\cN,\angle\eta_2),\ \Seg(\cN, P_{m}, \gamma_1)=0.\]

  For any consecutive pair $(a,b):=(i_k,i_{k+1})$ of $\hat{S}$, where $1\leq k\leq c$, we have
\[\Seg(\cM,P_{a-1}, \gamma_a)\neq 0\neq \Seg(\cM, P_{b-1},\gamma_b),
\]
\[\Seg(\cN,P_{a-1}, \gamma_a)\neq 0\neq \Seg(\cN, P_{b-1},\gamma_b),
\]
and
\[\Seg(\cM,P_{l-1}, \gamma_l)=0=\Seg(\cN,P_{l-1}, \gamma_l)
\]
 for $a<l<b$ (cf. Figure \ref{f:ab}). Moreover,
 \[\begin{cases}
 \Seg(\cM,\angle\eta_a)=d_{a-1}\\
 \Seg(\cM,\angle\eta_a)+\Seg(\cM,\angle\eta_{a+1})+\Seg(\cM,P_{a-1},\gamma_a)=d_{a}\\
 \Seg(\cM,\angle\eta_{a+1})+\Seg(\cM,\angle\eta_{a+2})=d_{a+1}\\
 \ \ \ \vdots\\
 \Seg(\cM,\angle\eta_{b-3})+\Seg(\cM,\angle\eta_{b-2})=d_{b-3}\\
 \Seg(\cM,\angle\eta_{b-2})=d_{b-2}\\
 \Seg(\cM,\angle\eta_{b})=d_{b-1}
 \end{cases}
 \]
 and 
  \[\begin{cases}
 \Seg(\cN,\angle\eta_a)=d_{a-1}\\
 \Seg(\cN,\angle\eta_a)+\Seg(\cN,\angle\eta_{a+1})+\Seg(\cN,P_{a-1},\gamma_a)=d_{a}\\
 \Seg(\cN,\angle\eta_{a+1})+\Seg(\cN,\angle\eta_{a+2})=d_{a+1}\\
 \ \ \ \vdots\\
 \Seg(\cN,\angle\eta_{b-3})+\Seg(\cN,\angle\eta_{b-2})=d_{b-3}\\
 \Seg(\cN,\angle\eta_{b-2})=d_{b-2}\\
 \Seg(\cN,\angle\eta_{b})=d_{b-1}
 \end{cases}.
 \]
 It follows that  $\Seg(\cM,\eta_i)=\Seg(\cN,\eta_i)$ for $a\leq i\leq b$ and $\Seg(\cM,P_{a-1},\gamma_a)=\Seg(\cN,P_{a-1},\gamma_a)$.
 By the arbitrariness of the pair $(i_k,i_{k+1})$, we conclude the desired equalities. This finishes the proof.  
\end{proof}
\begin{remark}
According to the proof of Lemma \ref{l:m-gons}, the conclusion of Lemma \ref{l:m-gons} holds for any $m(\geq 3)$-gon $\Delta$ provided that there is an arc $\gamma_i$ such that $\Seg(\cM, P_{i-1}, \gamma_i)\neq 0$.
\end{remark}

Now Theorem \ref{t:main-result-1} is a direct consequence of Lemma \ref{l: one loop} and Lemma \ref{l:m-gons}.

\section{Application to gentle algebras}\label{s:dimension-vector}

 \subsection{Gentle algebras }
 Let $Q$ be a quiver. Denote by $Q_0=\{1,\dots,n\}$  the set of vertices and $Q_1$ the set of arrows.  For an arrow $\alpha\in Q_1$, we denote by $s(\alpha)$ and $t(\alpha)\in Q_0$ the source and target of $\alpha$ respectively.

A finite-dimensional $k$-algebra $A$ is {\it gentle} if it admits a presentation $A=kQ/I$ satisfying the following conditions:
\begin{itemize} 
\item[(G1)] Each vertex of $Q$ is the source of at most two arrows and the target of at most two arrows;
\item[(G2)] For each arrow $\alpha$, there is at most one arrow $\beta$ (resp. $\gamma$) such that $\beta\alpha\in I$ (resp. $\gamma\alpha\not\in I$);
\item[(G3)] For each arrow $\alpha$, there is at most one arrow $\beta$ (resp. $\gamma$) such that $\alpha\beta\in I$ (resp. $\alpha\gamma\not\in I$);
\item[(G4)] I is generated by paths of length 2.
\end{itemize}

Let $A=kQ/I$ be a gentle algebra.
An oriented cycle $c=\alpha_s\alpha_{s-1}\cdots \alpha_1$ of length $s$ in the quiver $Q$ is said to have {\it full relations} if $\alpha_{i+1}\alpha_i\in I$ for $i=1,\dots, s-1$ and $\alpha_1\alpha_s\in I$. In particular, any loop is an oriented cycle with full relations. 

Let $e_i$ be the primitive idempotent of $A$ associated with the vertex $i\in Q_0$ and $P_i=e_iA$ the corresponding right projective $A$-module.
The {\it Cartan matrix} $C_A$ of $A$ has $\dimv P_i$ as its $i$-th column vector, \ie,
\[C_A=(\dimv P_1,\dots,\dimv P_n).
\]
The determinant of the Cartan matrix for a gentle algebra has been calculated by Holm \cite{H} explicitly. In particular, we have
\begin{theorem}\cite[Theorem 1]{H}\label{t:cartan-deter}
Let $A$ be a gentle algebra. The determinant of the Cartan matrix $C_A$ is zero if and only if $A$ admits an oriented cycle of even length with full relations.
\end{theorem}

\subsection{Geometric models of $\tau$-tilting theory of gentle algebras}

\begin{figure}[h]
\begin{tikzpicture}

\draw[thick,fill=black!20] (2,0.17)arc (250:290:3);
\node at (3,0.2){$p_\alpha=p_\beta$};
\fill(3,0) circle(1.5pt);
\draw[thick,color=blue] (3,0) parabola (1,-2);
\draw[thick,color=blue,dashed] (1,-2) arc (150:260:1);
\draw[thick,color=blue,->] (2.5,-1)--(2,-0.5);
\node at (2.7,-1.2){\tiny{$s(\alpha)$}};
\node at (1.9,-0.2){\tiny{$t(\alpha)$}};
\node at (2,-1) {$\alpha$};
\draw[thick,color=blue] (3,0) .. controls (2.6,-0.9) and (2.4,-1.1) .. (2,-2.2);
\draw[thick,color=blue] (3,0) .. controls (3.4,-0.9) and (3.6,-1.1) .. (4,-2.2);
\draw[thick,color=blue] (3,0) parabola (5,-2);
\draw[thick,color=blue,dashed] (5,-2) arc (30:-70:1);
\draw[thick,color=blue,->] (4,-0.5)--(3.5,-1);
\node at (4,-1) {$\beta$};
\node at (3.3,-1.2){\tiny{$t(\beta)$}};
\node at (4.1,-0.3){\tiny{$s(\beta)$}};

\draw[thick,fill=black!20] (9,0.17)arc (250:290:3);
\node at (10,0.2){$p_\alpha=p_\beta$};
\fill(10,0) circle(1.5pt);
\draw[thick,color=blue] (10,0) parabola (8,-2);
\draw[thick,color=blue,dashed] (9,-2.2) arc (160:220:1);
\draw[thick,color=blue,->] (9.5,-1)--(9,-0.5);
\node at (9,-1) {$\beta$};
\node at (9.7,-1.2){\tiny{$s(\beta)$}};
\node at (8.9,-0.2){\tiny{$t(\beta)$}};
\draw[thick,color=blue] (10,0) .. controls (9.6,-0.9) and (9.4,-1.1) .. (9,-2.2);
\draw[thick,color=blue] (10,0) .. controls (10.4,-0.9) and (10.6,-1.1) .. (11,-2.2);
\draw[thick,color=blue] (10,0) parabola (12,-2);
\draw[thick,color=blue,dashed] (11,-2.2) arc (20:-40:1);
\draw[thick,color=blue,->] (11,-0.5)--(10.5,-1);
\node at (11,-1) {$\alpha$};
\node at (10.3,-1.2){\tiny{$t(\alpha)$}};
\node at (11.1,-0.3){\tiny{$s(\alpha)$}};

\end{tikzpicture}
\caption{ Case $\alpha\beta=0$ when $p_\alpha=p_\beta$}\label{f:tiling-algebra}
\end{figure}
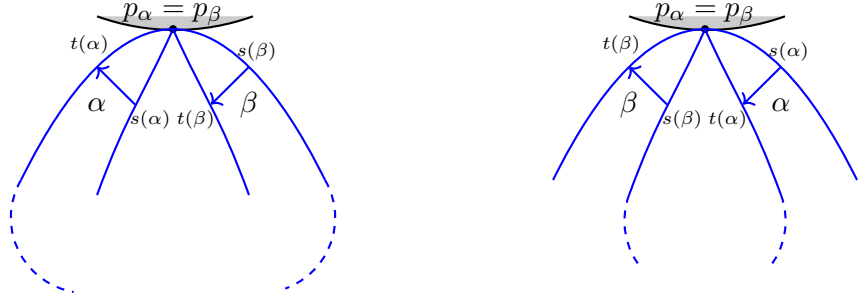

Let $(\mathbf{S}, \mathbf{M}, \mathbf{T})$ be a tiling.  Baur and  Sim\~{o}es \cite{BS} introduced a finite dimensional algebra $A^{\fT}$ called the {\it tiling algebra} associated with $(\mathbf{S}, \mathbf{M}, \mathbf{T})$. By construction, the tiling algebra $A^\fT$
is the factor algebra $kQ^{\mathbf{T}}/\langle R^{\mathbf{T}}\rangle$ of $kQ^{\mathbf{T}}$, where $(Q^{\mathbf{T}}, R^{\mathbf{T}})$ are described as follows.
\begin{itemize}
    \item The vertices in $Q^{\mathbf{T}}_0$ are indexed by the arcs in $\mathbf{T}$;
    \item There is an arrow $\alpha:i\to j\in Q_1^{\mathbf{T}}$ whenever the corresponding arcs $i$ and $j$ share an endpoint $p_\alpha\in \mathbf{M}$ such that $j$ follows $i$ anti-clockwise immediately;
    \item The relation set $R^{\mathbf{T}}$ consists of
    \begin{itemize}
        \item $\epsilon^2$, where $\epsilon$ is a loop;
        \item $\beta\alpha$ if $p_\alpha\neq p_\beta$, or the endpoints of the arc corresponding to $s(\beta)=t(\alpha)$ coincide and we are in one of the situations in Figure \ref{f:tiling-algebra}.
    \end{itemize}
\end{itemize}

The following result provides a geometric model for the representation theory of gentle algebras.
\begin{theorem}\cite[Theorem  2.10]{BS}\label{t:gentle-tiling}
Let $(\mathbf{S}, \mathbf{M},\mathbf{T})$ be a tiling. Then the tiling algebra $A^{\mathbf{T}}$ is a gentle algebra. Conversely, for any gentle algebra $A$, there is a tiling $(\mathbf{S}, \mathbf{M},\mathbf{T})$ such that $A\cong A^{\mathbf{T}}$ and $\mathbf{T}$ divides $\mathbf{S}$ into a collection of tiles of types $(I)$-$(V)$.
\end{theorem}

Let $A$ be a gentle algebra and $\tau$ the Aulander-Reiten translation. Recall that a finitely generated $A$-module $M$ is {\it $\tau$-rigid} if $\Hom(M,\tau M)=0$ (cf. \cite{AIR}). It is clear that every finitely generated projective $A$-module is $\tau$-rigid.
A pair $(M,P)$ is a {\it $\tau$-rigid pair} if $M$ is $\tau$-rigid and $\Hom_A(P,M)=0$, where $P$ is a finitely generated projective $A$-module.
Denote by $\opname{\tau-}\opname{rigid} A$ the set of isomorphism classes of $\tau$-rigid modules in $\mod A$, and by $\opname{ind\tau-}\opname{rigid} A$ the subset of $\opname{\tau-}\opname{rigid} A$ consisting of indecomposable ones.

According to Theorem \ref{t:gentle-tiling}, there is a tiling $(\mathbf{S}, \mathbf{M},\mathbf{T})$ such that $A\cong A^\fT$.
Each permissible arc segment $\eta$ satisfying condition (P2) gives an arrow $\alpha$ in  $Q^\fT$ or its formal inverse $\alpha^{-1}$. Furthermore, each permissible arc $\gamma$ gives a non-zero string and hence a string module of $A^{\fT}$. Clearly, the resulting string module depends only on the equivalence class of $\gamma$. We denote by $M([\gamma])$ the string $A^{\mathbf{T}}$-modules associated with the equivalence class $[\gamma]$ of $\gamma$.

The following correspondence is a consequence of \cite[Theorems 2.10 and 3.8]{BS} and  \cite[Propositions 5.3 and 5.6]{HZZ} (cf. also \cite[Lemma 2.4]{FGLZ}).
	
	\begin{lemma}\label{l:intersection-hom}
		Let $A=k Q/I$ be a finite-dimensional gentle algebra. Then there exists a tiling $(\fS,\fM,\fT)$ such that 
		\begin{itemize}
			\item we have a complete set of primitive orthogonal idempotents $\{e_{\ba}\mid \ba\in\fT\}$ of $A$ indexed by $\fT$, and
			\item there is a bijection
			\begin{eqnarray*}
			M\colon &\{\gamma\in\bP(\fT)\mid\Int(\gamma,\gamma)=0\}&\to\opname{ind\tau-}\opname{rigid} A,\\
			&\gamma&\mapsto M([\gamma])
			\end{eqnarray*}
			satisfying $\Intv_\fT(\gamma)=\dimv M([\gamma])$ for any  $\gamma\in\bP(\fT)$.
		\end{itemize}  Moreover, the bijection $M$ induces a bijection
		\[\mathscr{R}(\fT)\to\opname{\tau-}\opname{rigid} A,\]
		mapping $\mathcal{R}\in \mathscr{R}(\fT)$ to $\bigoplus\limits_{\gamma\in \mathcal{R}}M([\gamma])$  satisfying $\Intv_\fT(\mathcal{R})=\dimv \bigoplus\limits_{\gamma\in \mathcal{R}}M([\gamma])$.
	\end{lemma}

\subsection{Proof of Theorem \ref{t:main-result-2}}

``$(1)\Rightarrow (2)$": Suppose that $A$ admits an oriented cycle of even length with full relations. By Theorem \ref{t:cartan-deter}, the determinant of the Cartan matrix of $A$ is zero. Therefore there exist two non-isomorphic projective $A$-modules, say $P$ and $Q$, such that $\dimv P=\dimv Q$. Since $P$ and $Q$ are $\tau$-rigid, we obtain a  contradiction to the statement of $(1)$.

``$(2)\Rightarrow (1)$":
Let $(\fS, \fM,\fT)$ be a tiling such that $A\cong A^\fT$.
Since $A$ has no oriented cycle of even length with full relations, we conclude that  $(\fS,\fM,\fT)$ has no tiles of type $(II)$ nor  $m$-gon of type $(V)$ with $m$ even. Let $X,Y$ be two $\tau$-rigid $A$-modules such that $\dimv X=\dimv Y$. By Lemma~\ref{l:intersection-hom}, there are two multisets $\cM$ and $\cN$ such that $\cM=M^{-1}(X),\cN=M^{-1}(Y)$. Moreover,  $\Intv_{\fT}(\cM)=\dimv X=\dimv Y=\Intv_{\fT}(\cN)$. According to  Theorem~\ref{t:main-result-1}, we have $\cM=\cN$. By Lemma~\ref{l:intersection-hom} again, we conclude that $X=M(\cM)=M(\cN)= Y$. This completes the proof.

\section{Application to cluster algebras}\label{s:app-cluster}

\subsection{Recollection on cluster algebras}
Since  denominator vectors do not depends on the choice of coefficients, we only recall basic definitions and results for cluster algebras with trivial coefficients. We follow \cite{FZ07}. For an integer $b$, we use the notation $[b]_+=\max(b,0)$. 

Let $n$ be a positive integer and $\mathcal{F}$ the field of rational functions in $n$ variables. An $n\times n$ integer matrix $B$ is {\it skew-symmetrizable} if there is a diagonal matrix $S=\opname{diag}\{s_1, \dots, s_n\}$ with positive integers $s_1, \dots, s_n$ such that $SB$ is skew-symmetric. In this case, we call $S$ a {\it skew-symmetrizer} of $B$. A {\it labeled seed} in $\mathcal{F}$ is a pair $(B, \mathbf{x})$, where
\begin{itemize}
\item $B=(b_{ij})$ is an $n\times n$ skew-symmetrizable integer matrix, and
\item $\mathbf{x}=(x_1, \dots, x_n)$ is an $n$-tuple of elements of $\mathcal{F}$ forming a free generating set of $\mathcal{F}$.
\end{itemize}
We say that $\mathbf{x}$ is a {\it cluster} and refer to $x_i$ and $B$ as the {\it cluster variables} and the {\it exchange matrix} of the seed $(B, \mathbf{x})$, respectively.

Let $(B,\mathbf{x})$ be a labeled seed, and let $k \in\{1,\dots, n\}$. The \emph{seed mutation $\mu_k$ in direction $k$} transforms  $(B,\mathbf{x})$ into another labeled seed $\mu_k(B,\mathbf{x})=(B',\mathbf{x'})$ defined as follows:
\begin{itemize}
\item The entries of $B'=(b'_{ij})$ are given by 
\[
b'_{ij}=\begin{cases}-b_{ij} &\text{if $i=k$ or $j=k$,} \\ 
b_{ij}+\left[ b_{ik}\right] _{+}b_{kj}+b_{ik}\left[ -b_{kj}\right]_+ &\text{otherwise.}
\end{cases}
\]

\item The cluster variables $\mathbf{x'}=(x'_1, \dots, x'_n)$ are given by
\[
x'_j=\begin{cases}\dfrac{\mathop{\prod}\limits_{i=1}^{n} x_i^{[b_{ik}]_+}+\mathop{\prod}\limits_{i=1}^{n} x_i^{[-b_{ik}]_+}}{x_k} &\text{if $j=k$,}\\
x_j &\text{otherwise.}
\end{cases}
\]
\end{itemize}

Let $\mathbb{T}_n$ be the \emph{$n$-regular tree} whose edges are labeled by the numbers $1, \dots, n$ such that the $n$ edges emanating from each vertex have different labels. We write 
\begin{xy}(0,1)*+{t}="A",(10,1)*+{t'}="B",\ar@{-}^k"A";"B" \end{xy} 
to indicate that vertices $t,t'\in \mathbb{T}_n$ are joined by an edge labeled by $k$. 

 A {\it seed pattern} $t\mapsto \Sigma_t$ is an assignment of a seed $\Sigma_t=(B_t, \mathbf{x}_t)$ to each vertex $t$ of $\mathbb{T}_n$ such that the seeds
assigned to vertices $t$ and $t'$ linked by an edge labeled $k$ are obtained from each other
by the seed mutation $\mu_k$. For a given seed $(B, \mathbf{x})$, a seed pattern is uniquely determined by an assignment
of $(B, \mathbf{x})$ to a root vertex $t_0\in \mathbb{T}_n$.

Fix a seed pattern determined by assigning the seed $(B, \mathbf{x})$ to the root vertex $t_0$. Let $\Sigma_t=(B_t, \mathbf{x}_t)$ be the seed associated with $t\in \mathbb{T}_n$. The elements of $\Sigma_t$ are denoted as follows:
\[
\mathbf{x}_t=(x_{1;t},\dots,x_{n;t}),\ B_t=(b_{ij;t}).
\]
In particular, at $t_0$, we denote
\[
\mathbf{x}=\mathbf{x}_{t_0}=(x_1,\dots,x_n),\ B=B_{t_0}=(b_{ij}).
\]
A {\it cluster monomial} is a monomial in cluster variables all of which belong to the same cluster. 
The {\it cluster algebra $\cA(B)=\cA(B,\mathbf{x})$} associated with the seed pattern $t\mapsto \Sigma_t$ of $(B, \mathbf{x})$ is
the $\mathbb{Z}$-subalgebra of $\mathcal{F}$ generated by all the cluster variables $\mathcal{X}=\{x_{i;t}\}_{1\leq i\leq n;t\in \mathbb{T}_n}$. 
We refer to the cluster variables $x_1,\dots, x_n$ at the root vertex as the {\it initial cluster variables} and the seed $\Sigma_{t_0}$ as the {\it initial seed} of the cluster algebra $\mathcal{A}(B)$.

A cluster algebra $\cA$ is of {\it finite type} if $\cA$ has finitely many cluster variables. According to \cite{FZ03}, cluster algebras of finite type are classified by finite root system and are independent of the choice of coefficients.
Let $B=(b_{ij})$ be an $n\times n$ integer skew-symmetrizable matrix. The {\it Cartan counterpart} $A(B)=(c_{ij})$ of $B$ is an $n\times n$ integer matrix such that $c_{ii}=2$ for $i=1,\dots,n$, and $c_{ij}=-|b_{ij}|$ for $i\neq j$. In particular,  $A(B)$ is a generalized Cartan matrix.
By \cite{FZ03}, the following is well-defined.
\begin{definition}\
Let $(B,\mathbf{x})$ be a labeled seed and $\mathcal{A}$ the  cluster algebra associated with a seed pattern of $(B,\mathbf{x})$.
We say that $\cA$ is of type $\ast$ if there is a vertex $t\in \mathbb{T}_n$ such that $A(B_t)$ is of type $\ast$, where $\ast\in \{\mathbb{A}, \mathbb{B}, \mathbb{C},\mathbb{D},\mathbb{E}_6,\mathbb{E}_7,\mathbb{E}_8,\mathbb{F}_4,\mathbb{G}_2\}$. In this case, we also say that the exchange matrix $B$ is of type $\ast$.
\end{definition}
\subsection{Integer vectors arising from cluster algebras}
Throughout this section, let $\cA$ be a cluster algebra of rank $n$ with a fixed seed pattern $t\mapsto \Sigma_t$ whose root vertex is $t_0$. In particular, for each vertex $t\in \mathbb{T}_n$, we have a cluster $\mathbf{x}_t=(x_{1;t},\dots, x_{n;t})$ and an exchange matrix $B_t=(b_{ij;t})$. For $t_0$, we have $\mathbf{x}_{t_0}=(x_1,\dots, x_n)$ and $B_{t_0}=B=(b_{ij})$.

According to Laurent phenomenon \cite[Theorem 3.5]{FZ07}, every cluster variable $x_{j;t}\in \cA$ can be uniquely written as  
\[
x_{j;t} = \frac{N_{j;t}(x_1, \dots, x_n)}{x_1^{d_{1j;t}} \cdots x_n^{d_{nj;t}}},\quad d_{kj;t}\in \Z,
\]
where $N_{j;t}(x_1, \dots, x_n)$ is a polynomial with coefficients in~$\Z$
which is not divisible by any initial cluster variable~$x_j\in\mathbf{x}_{t_0}$. 
The integer vector 
\[\fd_{j;t}:=\begin{bmatrix}
d_{1j;t},\dots, d_{nj;t}
\end{bmatrix}^{\opname{tr}}\in \mathbb{Z}^n
\]
is called the {\it denominator vector}  (d-vector for short) of $x_{j;t}$ with respect to the initial cluster $\mathbf{x}_{t_0}$. 
The matrix $D_t:=(\fd_{1;t},\dots,\fd_{n;t})$ is called the {\it $D$-matrix} at vertex $t$.
For a cluster monomial $\mathbf{m}=x_{i;t}^{k_1}\cdots x_{n;t}^{k_n}$ of $\mathbf{x}_t$, the {\it denominator vector} $\fd(\mathbf{m})$ of $\mathbf{m}$ with respect to $\mathbf{x}_{t_0}$ is defined to be 
\[\fd(\mathbf{m})=\sum_{j=1}^nk_j\fd_{j;t}.
\]

In the following,
we recall the definitions of $c$-vectors, $g$-vectors and $f$-vectors via recursions and refer to \cite{FZ07} for the original definitions.

The family {\it $C$-matrices} $C_t=(c_{ij;t}):=(\mathbf{c}_{1;t},\dots,
 \mathbf{c}_{n;t})$ for $t\in \mathbb{T}_n$ can be defined by the initial conditions
 \[\mathbf{c}_{j;t_0}=\mathbf{e}_j\ (j=1,\dots,n)
 \]
together with the recurrence relations
 \begin{align*}
	\mathbf{c}_{j;t'} =
	\begin{cases}
		-\mathbf{c}_{j;t} & \text{if $j= k$;} \\[.05in]
		\mathbf{c}_{j;t} + [b_{kj;t}]_+ \ \mathbf{c}_{k;t} +b_{kj;t} [-\mathbf{c}_{k;t}]_+
		& \text{if $j\neq k$}.
	\end{cases}
\end{align*} 
 for any \begin{xy}(0,1)*+{t}="A",(10,1)*+{t'}="B",\ar@{-}^k"A";"B" \end{xy} in $\mathbb{T}_n$. Here $\mathbf{e}_1,\cdots, \mathbf{e}_n$ are unit vectors in $\mathbb{Z}^n$ and for an integer matrix $C$, we write $[C]_+$ for the matrix
obtained from $C$ by applying the operation $c\mapsto [c]_+$  to all entries of $C$. The $n\times n$ integer matrix $C_t$
is called the \emph{$C$-matrix} at vertex $t$ with respect to the root vertex $t_0$ and its column vectors are {\it $c$-vectors}.

The family of {\it $G$-matrix} $G_t=(\mathbf{g}_{1;t},\dots,g_{n;t})$ for $t\in \mathbb{T}_n$ can be defined by the initial conditions
\[\mathbf{g}_{j;t_0}=\mathbf{e}_j\ (j=1,\dots,n)
\] 
together with the recurrence relations
\begin{align*}
\mathbf{g}_{j;{t'}}&=\begin{cases}
\mathbf{g}_{j;t} \ \ & \text{if } j\neq k;\\
-\mathbf{g}_{k;t}+\mathop\sum\limits_{i=1}^{n}[b_{ik;t}]_+\mathbf{g}_{i;t}-\mathop\sum\limits_{i=1}^{n}[c_{ik;t}]_+\mathbf{b}_j &\text{if }  j=k.
\end{cases}
\end{align*}
 for any \begin{xy}(0,1)*+{t}="A",(10,1)*+{t'}="B",\ar@{-}^k"A";"B" \end{xy} in $\mathbb{T}_n$, 
where $\mathbf{b}_j$ is the $j$-th column of $B$. The column vector $\mathbf{g}_{j;t}$ is called the {\it $g$-vector} of the cluster variable $x_{j;t}$ with respect to the initial cluster $\mathbf{x}_{t_0}$. For a cluster monomial $\mathbf{m}=x_{i;t}^{k_1}\cdots x_{n;t}^{k_n}$ of $\mathbf{x}_t$, the {\it $g$-vector} $\mathbf{g}(\mathbf{m})$ of $\mathbf{m}$ with respect to $\mathbf{x}_{t_0}$ is defined to be 
\[\mathbf{g}(\mathbf{m})=\sum_{j=1}^nk_j\mathbf{g}_{j;t}.
\]

The family of {\it $F$-matrix} $F_t=(\mathbf{f}_{1;t},\dots, \mathbf{f}_{n;t})$ for $t\in \mathbb{T}_n$ can be defined  
 by the initial conditions
 \[\mathbf{f}_{j;t_0}=\mathbf{0} \ (j=1,\dots, n)
 \]
 together with the recurrence relations
\begin{align*}\label{f-recursion}
\mathbf{f}_{j;{t'}}&=\begin{cases}
\mathbf{f}_{j;t} \ \ & \text{if } j\neq k;\\
-\mathbf{f}_{k;t}+\max \left([\mathbf{c}_{k;t}]_++\mathop{\sum}\limits_{i=1}^n[b_{ik;t}]_+\mathbf{f}_{i;t},\ [-\mathbf{c}_{k;t}]_++\mathop{\sum}\limits_{i=1}^n[-b_{ik;t}]_+\mathbf{f}_{i;t}\right)\ \ &\text{if }  j=k.
\end{cases}
\end{align*}
 for any \begin{xy}(0,1)*+{t}="A",(10,1)*+{t'}="B",\ar@{-}^k"A";"B" \end{xy} in $\mathbb{T}_n$. 
The column vector $\mathbf{f}_{j;t}$ is called the  {\it $f$-vector} of $x_{j;t}$ with respect to the initial cluster $\mathbf{x}_{t_0}$. For a cluster monomial $\mathbf{m}=x_{i;t}^{k_1}\cdots x_{n;t}^{k_n}$ of $\mathbf{x}_t$, the {\it $f$-vector} $\mathbf{f}(\mathbf{m})$ of $\mathbf{m}$ with respect to $\mathbf{x}_{t_0}$ is defined to be 
\[\mathbf{f}(\mathbf{m})=\sum_{j=1}^nk_j\mathbf{f}_{j;t}.
\]

The following theorem collects some important properties of these vectors.

\begin{theorem}\label{t:property}
	\begin{enumerate}
		\item Different cluster monomials have different $g$-vectors.
		\item For each vertex $t\in \mathbb{T}_n$, we have $G_t^{\opname{tr}}SC_t=S$, where $S$ is a skew-symmetrizer of $B$.
	
		\item If $x_{i;t}\not\in \mathbf{x}_{t_0}$, then the denominator vector $\fd(x_{i;t})$ is  non-negative.
		\item The $k$-th component of $\fd(x_{i;t})$ depends only on $x_{i;t}$ and $x_{k;t_0}$, not on the cluster containing $x_{k;t_0}$; it is zero if and only if there is a cluster containing both $x_{i;t}$ and $x_{k;t_0}$.
		\item The cluster variable $x_{i;t}\not\in \mathbf{x}_{t_0}$ if and only if $\mathbf{f}_{i;t}$ is non-negative.
		
		\item A cluster $\mathbf{x}_t$ contains $x_k$ if and only if entries of the $k$th row of $F_t$ are all $0$.
		\item Assume moreover that the cluster algebra $\mathcal{A}$ is of finite type. If $x_{j;t}\not\in \mathbf{x}_{t_0}$, then $\mathbf{d}_{j;t}=\mathbf{f}_{j;t}$.
	\end{enumerate}
\end{theorem}
\begin{proof}
	(1) was proved by \cite{GHKK}(for cluster algebras of finite type,  it also follows from \cite{Demonet}). (2) follows from \cite[Thoerem 1.2]{NZ12} and \cite{GHKK}. (3) and (4) were proved by \cite[Theorem 6.3]{CL}. (5) and (6) are proved by \cite[Theorem 3.3 (1) and (4)]{FGy} and (7) follows from \cite[Theorem 1.8]{Gyoda} and $(3)$.
\end{proof}
\begin{remark}
In the following, we only use Theorem \ref{t:property} for cluster algebras of type $\mathbb{B}$ and $\mathbb{C}$. The corresponding results can be deduced from the additive categorifications of cluster algebras of type $\mathbb{B}$ and $\mathbb{C}$ (\cf \cite{Demonet,FGL21a,FGL21b}).
\end{remark}

Denote by $B^{\vee}:=-B^{\opname{tr}}$ and let $(B^\vee, \mathbf{x}^\vee)$ a labeled seed in $\mathcal{F}$. By assigning $(B^\vee, \mathbf{x}^\vee)$ to the vertex $t_0$, we obtain a seed pattern $t\mapsto \Sigma_t^\vee=(B_t^\vee, \mathbf{x}^\vee)$ of $(B^\vee, \mathbf{x}^\vee)$.
 Denote by $\mathcal{A}^\vee$ the associated cluster algebra, which is called the {\it Langlands dual} of $\mathcal{A}$. 
 For each vertex $t\in \mathbb{T}_n$, we also have a cluster $\mathbf{x}_t^\vee=(x_{1;t}^\vee,\dots, x_{n;t}^\vee)$,  $D$-matrix $D_t^\vee$, $C$-matrix $C_t^\vee$, $G$-matrix $G_t^\vee$ and $F$-matrix $F_t^\vee$ of the cluster algebra $\mathcal{A}^\vee$.
\begin{lemma}\label{l:duality}
Let $S$ be a skew-symmetrizer of $B$.
	For each vertex $t\in \mathbb{T}_n$, we have
	\begin{eqnarray*}F_t=S^{-1}F_t^\vee S;\\
	C_t=S^{-1}C_t^\vee S.
	\end{eqnarray*}
	
\end{lemma}
\begin{proof}
	The first equality was proved by \cite[Lemma 4.16]{FGy}, while the second equality was proved by \cite[(2.7)]{NZ12}.
\end{proof}

\subsection{Cluster algebras associated to marked surfaces with triangulations}\label{ss:cluster-algebra-surface-type}

Let $(\fS,\fM)$ be a marked surface. Recall that $(\fS, \fM)$ has no punctures. Let $\fT$ be a triangulation of $(\fS,\fM)$. Fomin, Shapiro and Thurston \cite{FST} constructed a skew-symmetric integer matrix $B_{\fT}$ associated with $\fT$. Therefore we have a cluster algebra $\mathcal{A}(B_\fT)$ associated with the triangulation $\fT$. We refer to \cite{FST} for the precise definition.
The following is a special case of \cite[Theorem 6.1]{FT}.
\begin{theorem}\cite[Theorem 6.1]{FT}\label{t:FT-correspondence}
Let $(\fS,\fM)$ be a marked surface with a triangulation $\fT$ and $\mathcal{A}(B_\fT)$ the associated cluster algebra. There is a bijection between the union of $\mathbb{P}(\fT)$ with arcs in $\fT$ and the set of cluster variables of $\mathcal{A}(B_\fT)$. Moreover, under this bijection, the initial cluster variables correspond to arcs in $\fT$ and  two cluster variables belong to the same cluster if and only if the corresponding arcs are compatible.
\end{theorem}
In the following, for a permissible arc $\delta$ in $\mathbb{P}(\fT)$, we will denote by $x_\delta$ the corresponding non-initial cluster variable.
The following is a direct consequence of \cite[Theorem 3.4]{ZZZ}, \cite[Theorem 1.5]{BZ} and \cite[Proposition 6.6]{FK} (cf. also \cite[Theorem 1.8]{Y}).
\begin{proposition}\label{p:f-vector}
For each permissible arc $\delta$ in $\mathbb{P}(\fT)$, the $f$-vector $\mathbf{f}(x_\delta)$ of $x_\delta$ with respect to the initial cluster  coincides with the intersection vector $\Intv_\fT(\delta)$ of $\delta$ with respect to $\fT$. 
\end{proposition}

As we have mentioned in the introduction, in order to record  the information of initial cluster variables in $f$-vectors, we need to modify the definition of $f$-vectors and we denote it by $\bar{\mathbf{f}}$. First, we define the $f$-vector $\bar{\mathbf{f}}(x_{i;t})$ for a cluster variable $x_{i;t}$ by setting
\[\bar{\mathbf{f}}(x_{i;t})=\begin{cases}\mathbf{f}(x_{i;t})& \text{if $x_{i;t}\neq x_k$ for any $1\leq k\leq n$ }\\
-\mathbf{e}_k& \text{if $x_{i;t}=x_k$ for some $k$}.
\end{cases}
\]
For a cluster monomial $\mathbf{m}=x_{1;t}^{k_1}\cdots x_{n;t}^{k_n}$, the {\it $f$-vector} $\bar{\mathbf{f}}(\mathbf{m})$ of $\mathbf{m}$ with respect to $\mathbf{x}_{t_0}$ is defined to be 
\[\bar{\mathbf{f}}(\mathbf{m})=\sum_{j=1}^nk_j\bar{\mathbf{f}}(x_{j;t}).
\]
According to Theorem \ref{t:property} (5) and (6), we can recover the $f$-vector $\mathbf{f}(\mathbf{m})$ from the one $\bar{\mathbf{f}}(\mathbf{m})$. The following is a reformulation of Theorem \ref{t:main-result-3}.
\begin{theorem}\label{t:reformulation-3}
Let $(\fS, \fM)$ be a marked surface with a triangulation $\fT$. Denote by $\mathcal{A}$ the associated cluster algebra. Let $\mathbf{m_1}$ and $\mathbf{m_2}$ be cluster monomials with $\bar{\mathbf{f}}(\mathbf{m}_1)=\bar{\mathbf{f}}(\mathbf{m}_2)$, then $\mathbf{m}_1=\mathbf{m}_2$. 
\end{theorem}
\begin{proof}
According to Theorem \ref{t:property} (5) and (6), it suffices to assume that $\mathbf{m}_1$ and $\mathbf{m}_2$ do not involve initial cluster variables. In this case, $\mathbf{f}(\mathbf{m}_i)=\bar{\mathbf{f}}(\mathbf{m}_i)$ for $i=1,2$. Now the result is a direct consequence of Theorem \ref{t:main-result-1}, Theorem \ref{t:FT-correspondence} and Proposition \ref{p:f-vector}.
\end{proof}

\subsection{Denominator conjecture for cluster algebras of type $\mathbb{A}$}
Let $\mathcal{A}$ be a cluster algebra of type $\mathbb{A}$ and $\mathbf{x}$ an initial cluster. It is well-known that there is a marked surface $(\fS,\fM)$ with a triangulation $\fT$ such that $\mathcal{A}(B_{\fT})\cong \mathcal{A}$. Moreover, $\fT$ corresponds to the initial cluster $\mathbf{x}$ under the bijection induced by Theorem \ref{t:FT-correspondence}.
On the other hand, by Theorem \ref{t:property} (7) and the definition of the $f$-vector $\bar{\mathbf{f}}$, we have $\bar{\mathbf{f}}(\mathbf{m})=\mathbf{d}(\mathbf{m})$ for any cluster monomial $\mathbf{m}$.
The following is a direct consequence of Theorem \ref{t:reformulation-3}.
\begin{corollary}\label{c:d-vector-A}
Let $\mathcal{A}$ be a cluster algebra of type $\mathbb{A}$. Different cluster monomials have different $d$-vectors with respect to any given cluster.
\end{corollary}

\subsection{Duality between cluster algebras of type $\bB$ and type $\bC$}
In this section, we show that the equivalence of denominator conjecture for cluster algebras of type $\bB$ and of type $\bC$, which reduces the proof of Theorem \ref{t:main-result-4} to the cases of type $\bA$ and $\bC$.

Throughout this section, let $B$ be  an $n\times n$  skew-symmetrizable integer matrix with a skew-symmetrizer $S=\opname{diag}\{s_1,\dots, s_n\}$. Let $(B,\mathbf{x})$ be a labeled seed in $\mathcal{F}$ and fix a seed pattern $t\mapsto \Sigma_t$ by assigning $(B,\mathbf{x})$ to the root vertex $t_0$. We also assign the labeled seed $(B^\vee, \mathbf{x}^\vee)$ to $t_0$ to obtain a seed pattern $t\mapsto \Sigma_t^\vee$ of $(B^\vee,\mathbf{x}^\vee)$.

The following weak version of denominator conjecture has been proved in \cite{FGL21a}.
\begin{lemma}\cite[Theorem 1.2]{FGL21a}\label{l:linearly-independeent}
	Assume that $B$ is of type $\mathbb{B}$ or $\mathbb{C}$, then for each vertex $t\in \mathbb{T}_n$, the column vectors of $D_t$ are linearly independent over $\mathbb{Q}$.
\end{lemma}

\begin{theorem}\label{t:duality-B-C}
	Denominator conjecture holds for cluster algebras of type $\mathbb{C}$ if and only if it holds for cluster algebras type $\mathbb{B}$.
\end{theorem}
\begin{proof}
	Let $B$ be of type $\mathbb{B}$, hence $B^\vee$ is of type $\mathbb{C}$.  We show that denominator conjecture holds for cluster algebra $\mathcal{A}(B)$ if and only if it holds for cluster algebra $\mathcal{A}(B^\vee)$.	
Without loss of generality,	we assume that $S=\operatorname{diag}\{1,2,\dots,2\}$ is a skew-symmetrizer of $B$. Consequently, $S_\bullet=\operatorname{diag}\{2,1,\dots,1\}$ is a skew-symmetrizer of $B^\vee$.

	We give a proof for the direction ``$\Leftarrow$" and the direction ``$\Rightarrow$" is similar.
	Assume that the denominator conjecture is true for the type $\mathbb{C}$ cluster algebra $\mathcal{A}(B^\vee)$, but it is not true for $\mathcal{A}(B)$. According to Lemma \ref{l:linearly-independeent}, there exists $t\neq s\in \mathbb{T}_n$ and two cluster monomials $u=x_{1;t}^{k_1}\cdots x_{n;t}^{k_n}$ and $v=x_{1;s}^{l_1}\cdots x_{n;s}^{l_n}$ such that 
\begin{itemize}
	\item $\fd(u)=\fd(v)$;
	\item $(k_1,\dots, k_n)>0$ and $(l_1,\dots, l_n)>0$;
	\item there exists $1\leq r\leq n$ such that $k_r>0$ and $x_{r;t}\not\in \mathbf{x}_s$.
\end{itemize}
By Theorem \ref{t:property} (3) and (4), we may assume moreover that both $\mathbf{x}_t$ and $\mathbf{x}_s$ contain no initial cluster variables.

Since $\fd(u)=\fd(v)$, we obtain
\[D_t[k_1,k_2,\dots, k_n]^{\opname{tr}}=D_s[l_1,l_2,\dots, l_n]^{\opname{tr}}.
\]
Applying Theorem \ref{t:property} (7),  we have
\[F_t[k_1,k_2,\dots, k_n]^{\opname{tr}}=F_s[l_1,l_2,\dots, l_n]^{\opname{tr}}.
\]
By Lemma \ref{l:duality}, 
\[S^{-1}F_t^\vee S[k_1,k_2,\dots, k_n]^{\opname{tr}}=S^{-1}F_s^\vee S[l_1,l_2,\dots, l_n]^{\opname{tr}}.
\]
Consequently, 
\begin{eqnarray}\label{f:F-matrix}
F_t^\vee [k_1,2k_2,\dots, 2k_n]^{\opname{tr}}=F_s^\vee [l_1,2l_2,\dots, 2l_n]^{\opname{tr}}.
\end{eqnarray}

Let $u^\vee=(x_{1;t}^\vee)^{k_1}(x_{2;t}^\vee)^{2k_2}\cdots(x_{n;t}^\vee)^{2k_n}$ and $v^\vee=(x_{1;s}^\vee)^{l_1}(x_{2;s}^\vee)^{2l_2}\cdots(x_{n;s}^\vee)^{2l_n}$ be cluster monomials of $\mathcal{A}(B^\vee)$. Recall that we have assume that $\mathbf{x}_t$ and $\mathbf{x}_s$ contains no initial cluster variables. By Theorem \ref{t:property} (5), each column vector of $F_t$ and $F_s$ is non-zero. Consequently, each column vector of $F_t^\vee$ and $F_s^\vee$ is non-zero. Hence both clusters $\mathbf{x}_t^\vee$ and $\mathbf{x}_s^\vee$ contains no initial cluster variables of $\mathcal{A}(B^\vee)$.
As a consequence, equality (\ref{f:F-matrix}) implies that $\fd(u^\vee)=\fd(v^\vee)$. It remains to show that $u^\vee$ and $v^\vee$ are different cluster monomials of $\mathcal{A}(B^\vee)$.

According to Theorem \ref{t:property} (2) and Lemma \ref{l:duality} (2), we obtain
\[ G_t^\vee=S_\bullet^{-1} G_t S_\bullet,
\]
\[G_s^\vee=S_\bullet^{-1} G_s S_\bullet.
\]
Let $G_t=(g_{ij;t})$ and $G_s=(g_{ij;s})\in \text{M}_n(\Z)$. If follows that 
\[G_t^\vee=\begin{bmatrix}
g_{11;t} &\frac{g_{12;t}}{2}&\cdots&\frac{g_{1n;t}}{2}\\ 2g_{21;t}&g_{22;t}&\cdots&g_{2n;t}\\
\vdots&\vdots&\vdots&\vdots\\
2g_{n1;t}&g_{n2;t}&\cdots&g_{nn;t}
\end{bmatrix}~\text{and}~G_s^\vee=\begin{bmatrix}
g_{11;s} &\frac{g_{12;s}}{2}&\cdots&\frac{g_{1n;s}}{2}\\ 2g_{21;s}&g_{22;s}&\cdots&g_{2n;s}\\
\vdots&\vdots&\vdots&\vdots\\
2g_{n1;s}&g_{n2;s}&\cdots&g_{nn;s}
\end{bmatrix}.
\]
We claim that $x_{r;t}^\vee\not\in \mathbf{x}_s^\vee$ and hence $u^\vee$ and $v^\vee$ are different cluster monomials.

If $r=1$, the g-vector $g_{1;t}^\vee$ of $x_{1;t}^\vee$ is $[g_{11;t},2g_{21;t},\cdots, 2g_{n1;t}]^{tr}$. If $x_{1;t}^\vee\in \mathrm{x}_s^\vee$, then there exists $1\leq j\leq n$ such that $g_{1;t}^\vee=g_{j;s}^\vee$. If $j=1$, then $g_{1;t}=g_{1;s}$, which implies that $x_{1;t}=x_{1;s}$ by Theorem \ref{t:property} (1), a contradiction. If $j>1$, then $g_{j;s}=2g_{1;t}$, which implies $x_{j;s}=x_{1;t}^2$, a contradiction.

One can discuss the case for $r>1$ similarly.
\end{proof}

\subsection{Denominator conjecture for cluster algebras of type $\bC$}\label{ss:proof-4}

Let $\cA$ be a cluster algebra of type $\bC$ of rank $n$ with a fixed seed pattern $t\mapsto \Sigma_t$. Let $t_0$ be the root vertex with initial seed $\Sigma_{t_0}=(B, \mathbf{x})$, where $B=(b_{ij})$. Without loss of generality, we assume that $S=\operatorname{diag}\{2,1,\dots,1\}$ is a skew-symmetrizer of $B$. 
We define a quiver $Q_B$ associated with $B$ as follows:
\begin{itemize}
\item the vertex set of $Q_B$ is $\{1,\dots, n\}$;
\item  there is a loop associated to the vertex $1$;
there are $b_{ij}$ arrows from vertex $i$ to vertex $j$ whenever $b_{ij}>0$ and $j\neq 1$; if $b_{i1}>0$, then we draw $\frac{b_{i1}}{2}$ arrows from vertex $i$ to vertex $1$.
\end{itemize}

The quiver $Q_B$ can  be characterized by the following conditions (cf. \cite{V,Yang}):
\begin{itemize}
\item[(a)] All non-trivial minimal cycles of length at least $2$ in the underlying graph is oriented and of length $3$;
\item[(b)] Any vertex has at most four neighbors;
\item[(c)] If a vertex has four neighbors, then two of its adjacent arrows belong to one $3$-cycle, and the other two belong to another $3$-cycle;
    \item[(d)]If a vertex has three neighbors, then two of its adjacent arrows belong to one $3$-cycle, and the third one does not belong to any $3$-cycle.
\item[(e)]There is a unique loop $\rho$ at the vertex $1$ which has one neighbor, or has two neighbors and its traversed by a $3$-cycle.
\end{itemize}
Let $A_B=kQ_B/I$ be the factor algebra of $kQ_B$ by the ideal $I$ generated by the square of the unique loop $\rho$ and all paths of length $2$ in a $3$-cycle. Clearly, $A_B$ is a gentle algebra.
\begin{proposition}\label{p:cat-C}
There is a bijection $\mathbb{X}$ between the set of isomorphism classes of indecomposable $\tau$-rigid $A_B$-modules and the set of non-initial cluster variables of $\mathcal{A}$. The bijection $\mathbb{X}$ induces a bijection between the set of $\tau$-rigid pairs of $A_B$ and the set of  cluster monomials of $\cA$. Moreover, for a $\tau$-rigid $A_B$-module $M$, denote by $\mathbb{X}_M$ the corresponding cluster monomial, we have
\[\fd(\mathbb{X}_M)=(\frac{m_1}{2},m_2,\dots, m_n)^{\text{tr}},
\]
where $\dimv M=(m_1,\dots, m_n)^{\text{tr}}$ is the dimension vector of $M$.
\end{proposition}
\begin{proof}
The bijection is a consequence of \cite{BMV} and \cite{CZZ, LX}. The equality of denominator vector was established by \cite[Theorem 1.3]{FGL21a}.
\end{proof}
\subsection*{Proof of Theorem \ref{t:main-result-4}}
According to Corollary \ref{c:d-vector-A} and Theorem \ref{t:duality-B-C}, it remains to prove the denominator conjecture for cluster algebras of type $\mathbb{C}$.
Let $\cA$ be a cluster algebra of type $\bC$ of rank $n$ with a fixed seed pattern $t\mapsto \Sigma_t$. Let $t_0$ be the root vertex with initial seed $\Sigma_{t_0}=(B, \mathbf{x})$.  Let $\mathbf{u}$ and $\mathbf{v}$ be two cluster monomials of $\cA$ such that $\fd(\mathbf{u})=\fd(\mathbf{v})$.  Recall that for the initial cluster variable $x_i$, we have $\fd(x_{i})=-\mathbf{e}_i$, $i=1,\dots,n$. According to Theorem \ref{t:property} (3)(4), it suffices to assume that both $\mathbf{u}$ and $\mathbf{v}$ do not admit any initial cluster variables in their expressions.

Let $A_B$ be the gentle algebra associated to $B$. By Proposition \ref{p:cat-C}, let $U$ and $V$ be the corresponding $\tau$-rigid $A_B$-modules of $\mathbf{u}$ and $\mathbf{v}$ respectively. By $\fd(\mathbb{X}_U)=\fd(\mathbb{X}_V)$, we obtain that $\dimv U=\dimv V$. On the other hand, according to the description of the quiver $Q_B$ of $A_B$ and its relations, we clearly know that $A_B$ does not have an oriented cycle of even length with full relations. It follows from Theorem \ref{t:main-result-2} that $U\cong V$. We conclude that $\mathbf{u}=\mathbf{v}$ by Proposition \ref{p:cat-C} .

\end{document}